\documentclass[11pt]{article}

\usepackage{amssymb,amsthm,amsmath,amsfonts}
\usepackage{graphicx}

\usepackage{subfigure}

\usepackage{lscape}
\usepackage{mathrsfs}

\usepackage{algorithmic}
\usepackage[ruled]{algorithm}

\usepackage{multirow}
\usepackage{verbatim}

\usepackage{color}

\newtheorem{lemma}{Lemma}
\newtheorem{theorem}{Theorem}

\def\no{\noindent}  \def\nb{\nonumber}
\def \Vh0{\stackrel{\circ}{V}_h}

\newcommand{\q}{\quad}    
\def\l{\label}

\def\ms{\medskip}

\def\bb{\begin{equation}} \def\ee{\end{equation}}

\usepackage{longtable}

\title{ Randomized Algorithms for Large-scale Inverse Problems \\ with General Regularizations}

\author{  Hua Xiang \thanks{ School of Mathematics and Statistics, Wuhan University, Wuhan 430072, P. R. China (hxiang@whu.edu.cn). }
\and Jun Zou \thanks{ Department of Mathematics,  The Chinese
University of Hong Kong, Shatin, New Territories, Hong Kong.
The work of this author
was substantially supported by Hong Kong RGC grants (Projects 405513 and 404611).
(zou@math.cuhk.edu.hk). } }

\begin{document}

\date{}
\maketitle

\begin{abstract}
We shall investigate randomized algorithms for solving large-scale
linear inverse problems with general regularizations. We first present some techniques
to transform inverse problems of general form into the ones of standard form,
then apply randomized algorithms to reduce large-scale systems of standard form
to much smaller-scale systems and seek their regularized solutions in combination with
some popular choice rules for regularization parameters.
%as we did in [H. Xiang, J. Zou, Regularization with
%randomized SVD for large-scale discrete inverse problems, Inverse
%Problem, 29 (2013) 085008].
Then we will propose a second approach to solve large-scale ill-posed systems
with general regularizations. This involves a new randomized generalized
SVD algorithm that can essentially reduce the size of the original large-scale ill-posed systems.
The reduced systems can provide approximate regularized solutions with
about the same accuracy as the ones by the classical generalized SVD, and more importantly,
the new approach gains obvious robustness, stability and computational time as it needs only to work on
problems of much smaller size.
Numerical results are given to demonstrated the efficiency of the algorithms.
\end{abstract}

%\begin{keywords}
%\end{keywords}

%\begin{AMS}
%\end{AMS}

%\tableofcontents

\section{Introduction}

Tikhonov regularization is one of the most popular and effective
techniques for the ill-conditioned linear system $K y = b$ arising
from the discretization of some linear or nonlinear inverse problems
\cite{BanksKunisch_Book89,ColtonKress_Book98,EnglHankeNeubauer_Book96},
where $K$  is an $m \times n$ matrix
%\footnote{Do all our results
%work for both $m>n$ and $m<n$ ? If not, we should make it very clear
%which case we will consider in this work. \textcolor{blue}{[[Our
%results work for both cases. The assumption of $m>n$ or $n<m$ never
%appears in the section 2 of standard form transformation. Only in
%section 2.4, we suppose that $m<n$ for RSVD. This is not essential.
%It works for $n<m$, but this RSVD is somewhat our own. And the case
%$m>n$ is already in literatures. ]]} }
and $b$ is an $m \times 1$
vector obtained from measurement data. The standard Tikhonov
regularization of this problem is of the form
\begin{equation} \label{ProbStandardForm}
\min_{y} || K y - b ||^2 + \mu^2 || y ||^2,
\end{equation}
where $\mu$ is the regularization parameter, and the 2-norm $\|\cdot\|$ is used
in this work unless otherwise specified. We will call formulation
\eqref{ProbStandardForm} as the standard form, where the second term is
for the regularization of the solution and the identity operator is used
for the regularisation. The identity regularisation is the simplest and most convenient one,
but it may not be the best in most applications. When we know
some additional a priori information about the physical
solution for a practical problem, we may apply some other more
effective regularisations. To differentiate from the standard one
\eqref{ProbStandardForm}, we will adopt other notation for
the general linear ill-conditioned system, namely $A x =b$, where $A$ is
an $m \times n$ matrix and $b$ is an $m \times 1$ vector. Then the
associated Tikhonov regularization in more general form is of the form
\begin{equation} \label{ProbGeneralForm}
\min_{x} ||A x - b||^2 + \mu^2 ||L x||^2,
\end{equation}
where the matrix $L$ is a $p \times n$ matrix, which may be a
discrete approximation to some differential operator, for  example,
the discrete Laplacian or gradient operator. When the null spaces of
$A$ and $L$ intersect trivially, i.e., $\mathcal{N}(A) \cap
\mathcal{N}(L) = \{0\}$, the regularised Tikhonov solution of
\eqref{ProbGeneralForm}  is unique.
As we shall see, the general form \eqref{ProbGeneralForm} can be transformed
into the standard one \eqref{ProbStandardForm} for very general regularisation
$L$.

%In some circumstances, we need to combine the above two
%regularization simultaneously to control not only the solution norm
%$||x||$, but also $||Lx||$. That is, we use a mixed regularization as follows.
%\begin{equation} \label{ProbMixedRegu}
%\min_{x} ||A x - b||^2 + \alpha^2 ||x|| + \beta^2 ||L x||^2.
%\end{equation}
%The usage of multiple penalties has proved very effective in
%promoting several distinct features [????].

When system \eqref{ProbStandardForm} or
\eqref{ProbGeneralForm} is large-scale, the traditional methods based on
singular value decomposition (SVD) or generalized SVD are very expensive and unstable,
and often infeasible for practical implementations.
For large-scale discrete systems \eqref{ProbStandardForm} of standard form, we can apply
the randomized SVD (RSVD) \cite{XiangZou_InverseProb13} to essentially reduce
the problem size, then combine the L-curve, GCV, and other discrepancy
principles to locate reasonable regularization parameters for solving the reduced regularisation
systems.
In this work, we shall focus on the solution of the general form
\eqref{ProbGeneralForm} of Tikhonov regularization. In section 2, we
discuss several techniques to transform the general form into the
standard one, then use the strategies in
\cite{XiangZou_InverseProb13} to solve the standard system. In
section 3, we consider the general form directly, and introduce a new
randomized generalized SVD (RGSVD) to reduce the problem size and then seek the regularised
solution.
Numerical experiments are given in section 4.

\section{Transformation into standard form and RSVD}

In this section we first discuss a general strategy to transform the problem
\eqref{ProbGeneralForm} of general form  into the standard one
\eqref{ProbStandardForm}, then apply the similar strategy as used in
\cite{XiangZou_InverseProb13}, which combines the randomized SVD with some
choice rules on regularization parameters, to solve the standard system.
%Most of the following standard form transformations
%are appeared in literatures, except the case 5.

\subsection{General transformation into the standard system}
\l{sec:transformation}
We demonstrate now how to transform the problem
\eqref{ProbGeneralForm} of general form with different regularisation
operators $L$ into the standard one \eqref{ProbStandardForm}.
Usually we assume that $\mathcal{N}(A) \cap \mathcal{N}(L) =\{0\}$,
such that the solution of \eqref{ProbGeneralForm} is unique. For the
cases where the matrix $L$ is of full column rank,
this assumption is automatically satisfied. But we shall consider
the most general case without this assumption, including both cases:
$\mathcal{N}(A) \cap \mathcal{N}(L) \neq
\{0\}$ and $\mathcal{N}(A) \cap \mathcal{N}(L) = \{0\}$. %And later consider some special cases of it.
For the case with $\mathcal{N}(A) \cap \mathcal{N}(L) \neq \{0\}$,
the solution of \eqref{ProbGeneralForm} is not unique, so
least-squares solutions with minimum norm will be sought.
In the sequel, we shall often use the Moore-Penrose generalized inverse $L^\dag$
of matrix $L$ \cite{GolubLoan_book13}.

We start with the following theorem which unifies the transformations
of problem \eqref{ProbGeneralForm} into the standard one \eqref{ProbStandardForm} for
all possibilities, and will discuss in section~\ref{sec:special_cases} several cases
when $A$ and $L$ have special structures,
where the results can be simplified.
% Please note that the solution in the following theorem has the same form as that in Theorem 1.
% But in the proof we need some other techniques, for example, the complete orthogonal decomposition \cite{GolubLoan_book13}.
%
\begin{theorem}\label{thm:general}
%\footnote{The theorem is exactly the same as Theorem 1 except for the opposite
%assumption about $L$. Why two opposite assumptions lead to the same
%results ? It is bit strange. If it is indeed the same, then we can
%simplify the statement of this theorem by simply citing the
%conditions inn Theorem 1. Please double check. \textcolor{blue}{I
%have already noticed that. The proof is not the same. And it is
%interesting that their appearances are similar. From another point
%of view, it is a good thing, so things can be unified.}
%
%Hua, since we have exactly the same results for two opposite cases:
%$\mathcal{N}(A) \cap \mathcal{N}(L) = \{0\}$ and $\mathcal{N}(A) \cap \mathcal{N}(L) \neq \{0\}$,
%so a natural question for readers will be that it is not necessary to discuss all other cases. In fact, since we have
%the same method and same results for two opposite cases, then this method works
%for all cases from the numerical point of view, and why do we need to discuss other cases ?
%The discussions are only necessary, unless more efficient methods can be applied for
%other individual cases. So here we need to have some convincing explanations to discuss why we have to discuss
%so many cases, why not simply to discuss only one case, and the result in Theorem 1 seems enough
%for all: only the technical proof is different but it does not affect the numerical method.}
%
% Suppose that $\mathcal{N}(A) \cap \mathcal{N}(L) \neq \{0\}$.
Let $W$, $Z$ and $P$ be any matrices satifying
$$
\mathcal{R} (W) = \mathcal{N}
(L), \q \mathcal{R} (Z) = \mathcal{R} (L),
\q Z^T Z=I,
\q
\mathcal{R}(P) = \mathcal{R}^\perp(AW), \q P^T P=I,
$$
and $L^\# = (I-W (AW)^\dag A) L^\dag$.
Then the least-squares solution with
minimum norm to the problem \eqref{ProbGeneralForm} of general form can be given by
\begin{equation} \label{Case6ThmSol}
 x_\mu = L^\# Z y_\mu + W (AW)^\dag b,
\end{equation}
where $y_\mu$ is the
minimizer of the following problem
\begin{equation} \label{Case6ThmOnesubprob}
\min_y  ||P^T A L^\dag Z y  - P^T b||^2 + \mu^2 || y ||^2.
\end{equation}
Equivalently,  $y_\mu$ can be obtained by solving 
%\footnote{This system looks not so convenient as it is minimising
%over the range space of $L$. I suggest to replace it by the system
%\eqref{Problem:StandardTransfromSubproblem}. In fact system
%\eqref{Problem:StandardTransfromSubproblem} is used in Algorithm 1,
%instead of system \eqref{Case6ThmOnesubprob2}.\textcolor{blue}{You
%are right. It's better to use (8). I have modified it here and the
%end of proof. Thanks!} }

\begin{equation} \label{Problem:StandardTransfromSubproblem}
\min_y \{ || A L^\# Z y -b ||^2 + \mu^2 || y ||^2 \}.
\end{equation}
%where $P$ satisfies $\mathcal{R}(P) = \mathcal{R}^\perp(AW)$ and $P^T P=I$.
\end{theorem}

\no {\it Proof}.
 Consider the SVD of matrix $L$:
\begin{equation}\label{eq:svdofL}
 L = U \Sigma V^T =
\begin{pmatrix} U_1, U_2 \end{pmatrix}
\begin{bmatrix} \Sigma_1 &
\\ & 0
\end{bmatrix}
\begin{pmatrix} V_1^T\\ V_2^T \end{pmatrix} =  U_1 \Sigma_1 V_1^T  ,
\end{equation}
 where $U=[U_1, U_2]$ and $V=[V_1, V_2]$ are
unitary matrices, $\Sigma_1$ includes the nonzero singular values in
the diagonal matrix $\Sigma$.
%To avoid introducing too many
%notations and follow the notation conventions, we abuse a little the
%notation $U$.

For any vector $x \in \mathbb{R}^n$, we write it as $x=V_1
\Sigma_1^{-1} \widehat{y} + V_2 \widehat{z}$, i.e., $x= L^\dag U_1
\widehat{y} + V_2 \widehat{z}$. Since
$span\{Z\}=\mathcal{R}(L)=span\{U_1\}$ and
$span\{W\}=\mathcal{N}(L)=span\{V_2\}$, we can set
$$x= L^\dag Z y +
W z .$$

%We still consider the SVD of $L$ as in \eqref{eq:svdofL} and seek the solution of the form $x=  L^\dag Z y + W z$.
%
%Since $\mathcal{N}(A) \cap \mathcal{N}(L) \neq \{0\}$, the matrix $AW$ is not of full rank.
Now we apply  the complete orthogonal decomposition on the
matrix $AW$ \cite{GolubLoan_book13}:
$$ Q^T (A W) \Pi = \begin{bmatrix} T_{11} & 0 \\ 0 & 0
\end{bmatrix}, $$
where $ Q$ and $\Pi$ are orthogonal matrices, and $T_{11}$ is
nonsingular matrix with the dimension determined by the rank of
$AW$.
For the case where $\mathcal{N}(A) \cap \mathcal{N}(L) = \{0\}$, the
matrix $AW$ is of full column rank, and the zero matrix on the right
side of $T_{11}$ will disappear. For the more general case with
$\mathcal{N}(A) \cap \mathcal{N}(L) \neq \{0\}$, the matrix $AW$ is
not of full column rank.
Let $Q=[Q_1, Q_2]$ be partitioned in compatible dimensions with
$T_{11}$. Then we have $span\{Q_1\} = \mathcal{R}(AW)$, and
$span\{Q_2\} =  \mathcal{R}^\perp(AW)$ % and $span\{Q_2\} = \mathcal{N}(AW)$
\cite{GolubLoan_book13}.
And it is easy to verify that
\begin{eqnarray*}
& & ||A x - b||^2 + \mu^2 ||L x||^2 \\
&=& ||Q^T A L^\dag  Z y + Q^T A W \Pi (\Pi^T z) - Q^T b||^2 + \mu^2 || L L^\dag Z  y ||^2 \\
&=& ||Q_1^T A L^\dag  Z y + T_{11} \bar{z}_1 - Q_1^T b||^2 + ||Q_2^T
A L^\dag Z y - Q_2^T b||^2 + \mu^2 ||  y ||^2 ,
\end{eqnarray*}
with $ \Pi^T z = [\bar{z}_1^T, \bar{z}_2^T ]^T$,
%\begin{bmatrix} \bar{z}_1 \\ \bar{z}_2 \end{bmatrix} $,
where $\bar{z}_1$ has the compatible dimension with $T_{11}$.
% and $Q=[Q_1, Q_2]$ is partitioned with the compatible dimension with $T_{11}$.
%
%The minimization of the formula above can be solved by the following two separate subproblems:

Minimizing the quadratic form above, we obtain
$\bar{z}_1 = T_{11}^{-1} Q_1^T (b - A L^\dag Z y_\mu )$, with $y_\mu$
given by the following subproblem of standard form
\begin{equation*}
\min_y  ||Q_2^T A L^\dag Z y - Q_2^T b||^2 + \mu^2 || y ||^2,
\end{equation*}
which is in fact the subproblem \eqref{Case6ThmOnesubprob}, since we
can see from the proof that $Q_2$ can be replaced by any orthonormal
matrix $P$ satisfying $span\{P\} = span\{Q_2\} =
\mathcal{R}^\perp(AW)$.

Now we can verify that the solution of problem \eqref{ProbGeneralForm} is given by
$$ x_\mu = L^\dag Z y_\mu + W \Pi \begin{bmatrix} \bar{z}_1 \\ \bar{z}_2 \end{bmatrix}  \quad \forall ~ \bar{z}_2.  $$
Noting $||x_\mu||^2=||L^\dag Z y_\mu||^2 +
\|[\bar{z}_1^T, \bar{z}_2^T ]^T\|^2$, we obtain the solution with minimum norm with $\bar{z}_2=0$:
\begin{eqnarray*}
 x_\mu &=& L^\dag Z y_\mu + W \Pi \begin{bmatrix}\bar{z}_1 \\ 0 \end{bmatrix}  \\
   & = & L^\dag Z y_\mu + W \Pi \begin{bmatrix} T_{11}^{-1} & 0 \\ 0 & 0
   \end{bmatrix} \begin{bmatrix} Q_1^T \\ Q_2^T \end{bmatrix} (b - A L^\dag Z y_\mu ) \\
   & = & L^\dag Z y_\mu + W (A W)^\dag (b - A L^\dag Z y_\mu ) \\
   & = & (I- W(A W)^\dag A) L^\dag Z y_\mu + W (A W)^\dag b \\
   & = & L^\# Z y_\mu + W (A W)^\dag b.
 \end{eqnarray*}
We can see that for the case with $\mathcal{N}(A) \cap
\mathcal{N}(L)=\{0\}$, the vector $\bar{z}_2$ disappears
automatically, so the above solution is the unique least-squares
solution of \eqref{ProbGeneralForm}.

We can easily verify that the minimizer of \eqref{Case6ThmOnesubprob} is given by
$$ 
y_\mu = [ Z^T (A L^\dag)^T P P^T (A L^\dag) Z + \mu^2 I ]^{-1} Z^T( A L^\dag)^T  P P^T b.
$$

Note that $PP^T$ is the orthogonal projection onto
$\mathcal{R}^\perp (AW)$. From the uniqueness of orthogonal
projection, we have $PP^T = I - AW (AW)^\dag$.
It is straightforward to check that
\begin{eqnarray*}
P P^T A L^\dag &=& [I - AW (AW)^\dag] A L^\dag = [A - AW (AW)^\dag A] L^\dag =  A L^\#\,,\\
(A L^\#)^T(A L^\# )
%&=& {L^\dag}^T (I - A^T {(AW)^\dag}^T W^T) A^T A (I - W (AW)^\dag A) L^\dag \\
%&=& {L^\dag}^T (A^T - A^T {(AW)^\dag}^T W^T A^T)  (A - A W (AW)^\dag A) L^\dag \\
%&=& {L^\dag}^T A^T(I -  {(AW)^\dag}^T (AW)^T)  (I - A W (AW)^\dag ) A L^\dag \\
%&=&(A L^\dag)^T (I -  AW (AW)^\dag)^T  (I - A W (AW)^\dag ) A L^\dag \\
%&=& (A L^\dag)^T (I - A W (AW)^\dag ) A L^\dag \\
 &=&  (A L^\dag)^T PP^T A L^\dag .
\end{eqnarray*}
Hence the solution can be rewritten as
\begin{equation*}
y_\mu =  [ Z^T(A L^\#)^T(A L^\#)Z + \mu^2 I ]^{-1} Z^T(A L^\#)^T b ,
\end{equation*}
which is obviously the minimizer of
\eqref{Problem:StandardTransfromSubproblem}.   
\qed

We know that the columns of matrix $Z$ span the range of $L$
and $Z^T Z=I$, any vector $\widehat{y} \in \mathcal{R}(L)$ can be
expressed as $\widehat{y} = Z y$. So the problem
\eqref{Problem:StandardTransfromSubproblem} is equivalent to
\begin{equation}\label{Case6ThmOnesubprob2}
\min_{\widehat{y} \in \mathcal{R}(L)} \{ || A L^\# \widehat{y} -b
||^2 + \mu^2 || \widehat{y} ||^2 \}\,,
\end{equation}
whose minimizer  is given by
\begin{equation*}
\widehat{y}_\mu = Z [ Z^T(A L^\#)^T(A L^\#)Z + \mu^2 I ]^{-1} Z^T(A
L^\#)^T b = Z  y_\mu .
\end{equation*}

%Note that we can also express $\bar{z}_1$ as $\bar{z}_1 =
%T_{11}^{-1} Q_1^T (b - A L^\dag \widehat{y} )$ with $\widehat{y}$ is given by
%\begin{equation*}
%\min_{\widehat{y} \in \mathcal{R}(L)}  ||Q_2^T A L^\dag \widehat{y}
%- Q_2^T b||^2 + \mu^2 ||\widehat{y} ||^2.
%\end{equation*}
%
%If we can find a matrix $W$ whose columns span the null space of
%$L$, that is, $span\{W\}=span\{V_2\}= \mathcal{N}(L)$, then we can
%use $W$ to replace $V_2$ in above derivation,

%Then the minimum norm least squares solution of general form problem
%\eqref{ProbGeneralForm} can be  rewritten as
%\begin{equation*}
% x = (I-W (AW)^\dag A) L^\dag \widehat{y} + W (AW)^\dag b.
%\end{equation*}
%***********************************************************************

\subsection{Practical realisation of the transformation}

%\footnote{Please add a remark here about the generation of matrices
%$W, Z$ and $P$ in Theorems 1-3. The readers may think they are
%expensive to generate so our transformations will be expensive.
%Mention some efficient ways to generate them. \textcolor{blue}{I add
%one page about this in the following}. }

%\footnote{It may be more clear to separate this subsection into
%subsections: the first one for the generations of matrices $W$, $Z$,
%and $P$ etc., the second one for RSVD. \textcolor{blue}{ I separate
%this part from RSVD, and put it as a subsection, followed by a new
%subsection on special cases. Do you think it is a proper way to
%treat this? Please note that the discussion of this subsection is
%somewhat general. For the special cases, there exists special
%treatments, but equivalent to the statements in Theorem 1. } }
Theorem~\ref{thm:general} gives a unified transformation that works for all
possible choices of regularisation matrix $L$ in \eqref{ProbGeneralForm}.
We now discuss some practical realisation of the matrices $W$, $Z$, $P$
and the oblique pseudoinverse $L^\#$ involved in the transformation as stated in Theorem~\ref{thm:general}.
By means of the standard SVD \eqref{eq:svdofL} of the matrix $L$,
we can choose $Z=U_1$ and
$W=V_2$ such that $\mathcal{R}(Z) = \mathcal{R}(L)$ and
$\mathcal{R}(W) = \mathcal{N}(L)$. But the SVD is rather expensive. Instead
we may use the complete orthogonal factorization \cite{GolubLoan_book13}
in practical computations when $L$ is not of full rank:
$$ U^T L V = \begin{bmatrix} T & 0 \\ 0 & 0
\end{bmatrix}, $$
where $U=[u_1, \cdots , u_p]$ and  $V=[v_1, \cdots , v_n]$ are
orthogonal matrices, and $T$ is a $r \times r$ nonsingular matrix,
with $r=rank(L)$.
Then we have \cite{GolubLoan_book13}
$$
\mathcal{R} (L) = span\{ u_1, \cdots, u_r\} , \q
\mathcal{N} (L) = span\{v_{r+1}, \cdots, v_n \}.
$$
So we can choose $W= [v_{r+1}, \cdots, v_n]$ and $Z=[u_1, \cdots,
u_r]$.
When matrix $L$ is of full rank,  the matrices $W$ and $Z$ can be determined by QR, or
QR with column pivoting, which are special cases of the complete
orthogonal factorization.

For the choice of matrix $P$, we perform QR with column pivoting on the matrix $AW$
$$ (AW) \Pi = [Q_1, Q_2] \begin{bmatrix} T_1 \\ 0 \end{bmatrix} =
Q_1 T_1 ,$$ where $\Pi$ is a permutation matrix, $T_1$ is of full
row rank, and $[Q_1, Q_2]$ is an orthogonal matrix.
%For the case
%where $\mathcal{N}(A) \cap \mathcal{N}(L) =\{0\}$, the matrix $T_1$
%is nonsingular upper triangular.
%
Then we have $span\{Q_1\} = \mathcal{R}(AW)$, and $span\{Q_2\} =
\mathcal{R}^\perp(AW)$. So we can choose $P=Q_2$ in Theorem~\ref{thm:general}.
On the other hand, we know from the proof of Theorem~\ref{thm:general} that problem
\eqref{Case6ThmOnesubprob} is the same as the minimisation \eqref{Problem:StandardTransfromSubproblem}.
%\begin{equation*}
% \min_y \{ ||
%A L^\# Z y -b ||^2 + \mu^2 || y ||^2 \}\,.
%\end{equation*}
So if we choose to solve system \eqref{Problem:StandardTransfromSubproblem} instead of \eqref{Case6ThmOnesubprob},
we can get rid of matrix $P$ in all computations.

For the oblique pseudoinverse $L^\#$, it involves the Moore-Penrose
inverse of $AW$ and can be computed as follows:
\begin{equation} \label{eq:computeAW}
(AW)^\dag = (Q_1 T_1 \Pi^T)^\dag  = \Pi (Q_1 T_1)^\dag = \Pi
T_1^\dag Q_1^T = \Pi T_1^T(T_1 T_1^T)^{-1} Q_1^T  .
\end{equation}
For the special case with $\mathcal{N}(A) \cap \mathcal{N}(L) =
\{0\}$, the matrix $A W$ is of full column rank, and we can use QR
factorization. Correspondingly  we have $\Pi=I$ and $(A W)^\dag =
T_1^{-1} Q_1^T$.
For most applications, the dimension of null space
$\mathcal{N}(L)$ is very low, for example, $\mathcal{N}(L)$ may be
spanned simply by a single vector $[1, 1, \cdots,1]^T $ or $[1, 2,
\cdots, n]^T$. So the matrix $AW$ is very tall skinny, $T_1$ is
a very small matrix, and the cost for computing the
Moore-Penrose inverse of  $AW$ is negligible.

Now we can summarise the solution to the problem
\eqref{ProbGeneralForm} of general form in Algorithm~1.
\begin{algorithm}\l{alg:general}
\caption{(Standard form transformation). }
\label{algo:StandardformTransformation}
\begin{algorithmic}
   \STATE{1. Generate matrix $W$ satisfying $\mathcal{R}(W)=\mathcal{N}(L)$.
   }
   \STATE{2. Generate matrix $Z$ satisfying $\mathcal{R}(Z)=\mathcal{R}(L)$ and $Z^T Z=I$.
   }
   \STATE{3. Compute $(AW)^\dag$, for example, by \eqref{eq:computeAW}.
   }
   \STATE{4. Form the generalized inverse $L^\# = (I-W (AW)^\dag A) L^\dag$.
   }
   \STATE{5. Find the minimizer $y_\mu$ of subproblem \eqref{Problem:StandardTransfromSubproblem}.
   }
   \STATE{6. Form the solution $x_\mu = L^\# Z y_\mu + W (AW)^\dag b$.
   }
\end{algorithmic}
\end{algorithm}

As we may see in the next section~\ref{sec:special_cases},
some steps of Algorithm
\ref{algo:StandardformTransformation} can be omitted for matrices
$L$ and $A$ of special properties, which are listed below:

\begin{enumerate}
\item When $L$ is of full column rank, we have $W=0$, hence Steps 1 and 3 can be dropped,
and the terms involving $W$ do not appear in Steps 4 and 6.

\item When $L$ is of full row rank, we have
$Z=I$ and can skip Step 2.

\item  If $\mathcal{N}(L) \subseteq \mathcal{N}(A)$,   we have $AW=0$,
hence Steps 1 and 3 can be dropped,
and the terms involving $AW$ do not appear in Steps 4 and 6.

\item When $L$ is of full row rank and $\mathcal{N}(A) \cap
\mathcal{N}(L) = \{0\}$, it is unnecessary to form the pseudo-inverses $(AW)^\dag$ and $L^\#$
explicitly, instead we can solve the subproblem \eqref{Case6ThmOnesubprob}.

\item If $L$ is rank-deficient and $\mathcal{N}(A) \cap
\mathcal{N}(L) = \{0\}$, then $AW$ is of full column rank, and the Moore-Penrose inverse $(AW)^\dag$ can be
directly achieved by a QR decomposition.

\item
One may use iterative methods to avoid forming the matrix $L^\#$ explicitly in Step 4.
Instead we need only to have a solver for the linear system $Lg=h$ with given right-hand sides
$h$ to achieve $g=L^\dag h$ approximately.
\end{enumerate}

For Step 5, one may apply the randomised SVD to first reduce the system size
essentially, then solve the reduced system in combination with some strategies
for regularization parameters, as did in \cite{XiangZou_InverseProb13}.
The standard form transformation described above is an effective and
efficient approach for solving the ill-posed problem
\eqref{ProbGeneralForm}, provided that the operations with $L^{-1}$,
$L^\dag$ or $L^\#$ can be efficiently implemented. When $L$ is the
discrete Laplacian, the actions of inverses can be done by the
algebraic multigrid method efficiently \cite{RugeStuben_AMG1987}.
For the cases where $L$ is of special structures, such as Hankel or Toeplitz,
there exist many fast solvers for
implementing the operations with $L^{-1}$ or $L^\dag$ \cite{chan07}.

%-------------------------------------------------------------------------

\subsection{Special cases}
\label{sec:special_cases}

In this subsection, we consider a few important special cases.
Although all these cases have the solutions of same form \eqref{Case6ThmSol}
(see Theorem~\ref{thm:general}), the solutions may be
realised very differently as it is shown below.

{\bf Case 1: $\mathcal{N}(A) \cap \mathcal{N}(L) \neq \{0\}$, and
$\mathcal{N}(L) \subseteq \mathcal{N}(A)$}.
This happens in certain practical problems, for example, the
lead-field matrix and the Laplacian have a vector of all ones in
their null spaces for the inverse problem from electrocardiography. 
For this case, $AW=0$
and $L^\# = L^\dag$. The least-squares solution with minimum norm is
given by
$x_\mu =L^\dag Z y_\mu$,  
where $y_\mu$ solves the minimisation of standard form
\begin{equation*}
\min_y  || A L^\dag Z y  - b||^2 + \mu^2 || y ||^2.
\end{equation*}

{\bf Case 2: $L$ is of full row rank, and $\mathcal{N}(A) \cap
\mathcal{N}(L) = \{0\}$}.
Different from the transformation used in Theorem~\ref{thm:general}, there is
an alternative approach in \cite{Elden_BIT77}, which
applies the following two QR decompositions:
$$ L^T = Q \begin{pmatrix} R \\ 0 \end{pmatrix}, \quad Q = [Q_1, W];
 \qquad   A W = U \begin{pmatrix} T \\ 0 \end{pmatrix}, \quad U = [U_1,
P], $$
where $Q$, $U$ are orthogonal matrices,  $R$ and $T$ are nonsingular
upper triangular matrices. For this case, the
matrix $Z$ can be chosen as the identity. So the three
matrices $W$, $P$ and $Z$ in Theorem~\ref{thm:general} are all well defined.
Next we shall derive the solution of \eqref{ProbGeneralForm}.
Note that the
range $\mathcal{R} (W) = \mathcal{R}^\perp(L^T) = \mathcal{N}(L)$,
hence $A W$ is of full rank, and $\mathcal{R}(P) = \mathcal{R}^\perp
(AW)$.
Suppose $x = L^\dag y + W z$, with $L^\dag = L^T (L L^T)^{-1} = Q_1
R^{-T}$, we can derive that
\begin{eqnarray*}
& & ||A x - b||^2 + \mu^2 ||L x||^2 \nb\\
&=& || A L^\dag y + A W z - b||^2 + \mu^2 || L L^\dag y ||^2 \nb\\
&=& ||U^T A L^\dag y + U^T A W z - U^T b||^2 + \mu^2 || y ||^2 \nb\\
&=& ||U_1^T A L^\dag y + T z - U_1^T b||^2 + ||P^T A L^\dag y - P^T
b||^2 + \mu^2 || y ||^2 . \l{eq:yz}
\end{eqnarray*}
Then we can show the minimisation  \eqref{ProbGeneralForm}  is
equivalent to the following two separated subproblems:
%\footnote{These two subproblems are not equivalent to the
%above coupled system for both $z$ and $y$, yes ? If it is not
%equivalent, then the following derivation may not make much sense.
%\textcolor{blue}{Answer: these two subproblems are equivalent to the
%above coupled system. The original coupled one can be separated. And
%the following derivation is meaningful. Generally, x=x(y,z),
%$\min_x$ cannot be separated by $\min_y + \min_z$. But for this
%special case, it is OK. }}
%Then $u$ can be derived by solving
\begin{equation*}
\min_y  ||P^T A L^\dag y - P^T b||^2 + \mu^2 || y ||^2 ~~\mbox{and}
~~ \min_z ||U_1^T A L^\dag y + T z - U_1^T b|| .
\end{equation*}
The first subproblem is the same as \eqref{Case6ThmOnesubprob} in
Theorem 1 corresponding to the matrix $Z=I$. We
can compute $z = T^{-1} U_1^T (b - A L^\dag y_\mu )$,
where $y_\mu$ is the minimizer of the first subproblem. Hence,
$$ x_\mu = L^\dag y_\mu + W z = L^\dag y_\mu + W T^{-1} U_1^T (b - A L^\dag
y_\mu).$$
Though the matrix $(AW)^\dag$ does not appear
explicitly, this solution is in fact equivalent to
\eqref{Case6ThmSol}, by the fact that $(AW)^\dag = T^{-1} U_1^T$ and  $Z=I$.
The solution of this case can be rewritten as
\begin{equation*} %\label{Tranform2Stand_case3}
 x_\mu = (I-W
(AW)^\dag A) L^\dag  y_\mu + W (AW)^\dag b = L^\# y_\mu + W
(AW)^\dag b.
\end{equation*}

For the QR factorization of large matrices, we may use the recently
developed new technique, communication-avoiding QR (CAQR), which invokes
tall skinny QR (TSQR) for each block column factorization, to speed
up the computation
\cite{DemmelGrigoriHoemmenLangou_SISC12,DemmelGrigoriGuXiang_SIMAX}.

{\bf Case 3: $L$ is of full column rank}. %, and $N(A) \cap N(L) = \{0\}$.
Using the skinny QR decomposition $L=Q_1 R$, where $R$ is
nonsingular and upper triangular,  and $Q_1$ is column orthogonal,  we
have $||L x|| = ||R x||$. Hence the problem \eqref{ProbGeneralForm}
of general form  is equivalent to the following system
\begin{equation} \label{case4_fullcolumnrank}
\min_x ||A R^{-1} R x - b||^2 + \mu^2 ||R x||^2.
\end{equation}
Then we can easily transform the system \eqref{ProbGeneralForm} to
the standard form \eqref{ProbStandardForm} by using $y=Rx$ and $K=A
R^{-1}$. This is efficient for practical computing since we need
only a skinny QR decomposition and a upper triangular solver. The
problem \eqref{case4_fullcolumnrank} is actually the same as the
problem \eqref{Case6ThmOnesubprob} in Theorem~\ref{thm:general}, by noting the facts
that $W=0$, $P=I$, $Z=Q_1$ and $L^\#=L^\dag = (L^TL)^{-1}L^T =  R^{-1}
Q_1^T$ in this case.

%In order to unify the subsequent results by using the generalized
%inverse of $L$, we can also use another transformation. Using the
%skinny QR decomposition $L=Q_1 R$, the general problem
%\eqref{ProbGeneralForm} is mathematically equivalent to
%\begin{equation*}
% \min_x ||A R^{-1} Q_1^T Q_1 R x - b||^2 + \mu^2 ||Q_1 R x||^2,
%\end{equation*}
%
%which can be transformed into the standard form
%\eqref{ProbStandardForm} simply by the means of the transformation
% $y =Q_1 R x$ and $K=A R^{-1} Q_1^T$.
%
%Noting the Moore-Penrose generalized inverse $L^\dag = (L^T L)^{-1}
%L^T = R^{-1} Q_1^T$, the transformation can be simplified to
%\begin{equation}\label{Tranform2Stand_case2}
% K=A L^\dag, \quad x= R^{-1} Q_1^T y = L^\dag y .
%\end{equation}

{\bf Case 4: $L$ is a square and nonsingular matrix}.
As $L$ is nonsingular, we can simply set
\begin{equation}\label{Tranform2Stand_case1}
K = A L^{-1},  \q x = L^{-1} y,
\end{equation}
then the problem \eqref{ProbGeneralForm} is rewritten in the
standard form \eqref{ProbStandardForm}.
% Then as in \cite{XiangZou_InverseProb13} we use the random strategies to
% solve this discrete ill-posed problem, together with L-curve
%or GCV for the regularization parameter determination.
%
The transformation \eqref{Tranform2Stand_case1} is applicable
whenever the actions of $L^{-1}$ can be performed efficiently. This
is the case when $L$ is sparse, banded, or of some special
structure.

\ms As we have seen from the above cases, by using the
transformation that may involve (generalized) matrix inverses,
%like the ones in \eqref{Tranform2Stand_case1}, \eqref{Tranform2Stand_case2}, \eqref{Tranform2Stand_case3}, \eqref{Case4_PCHansen_SolWrong1}, \eqref{Case4anotherderivSol1} or \eqref{Case5_lsqsol2},
we can transform the problem \eqref{ProbGeneralForm} of general form
into the problem \eqref{ProbStandardForm} of standard form. Then existing methods
for the standard form can be applied as we discuss in
the next subsection.

\subsection{Solution of the standard system \eqref{ProbStandardForm} by randomized SVD}

As we have seen in subsections~\ref{sec:transformation}-\ref{sec:special_cases},
the regularized solution $x_\mu$ of general form \eqref{ProbGeneralForm} can
be reduced to the solution to the standard system
\eqref{ProbStandardForm}. When the standard system \eqref{ProbStandardForm} is large-scale,
we can first apply randomized SVD algorithm (see Algorithm \ref{algo:RSVD*}
for $m<n$)
to reduce it to a much smaller system, then solve it by combining
with some existing choice rules for regularization parameters
\cite{XiangZou_InverseProb13}. Similar algorithm can be formulated
for $m>n$ \cite{HalkoMartinssonTropp_SIREV11}.

\begin{algorithm}
\caption{(RSVD). Given  $K \in \mathbb{R}^{m \times n} (m < n)$ and
$l<m$, compute an approximate rank-$l$ SVD:  $K \approx U \Sigma
V^T$  with $U \in \mathbb{R}^{m\times l}$, $\Sigma \in \mathbb{R}^{l
\times l}$ and $V \in \mathbb{R}^{n \times l}$. } \label{algo:RSVD*}
\begin{algorithmic}
   \STATE{1. Generate an $l \times m$ Gaussian random matrix $\Omega $ .
   }
   \STATE{2. Compute the $l \times n$ matrix $Y = \Omega K $.
   }
   \STATE{3. Compute the $n \times l$ orthonormal matrix $Q$ via QR factorization $Y^T = Q R $.
   }
   \STATE{4. Form the $m \times l$ matrix $B = K Q$.
   }
   \STATE{5. Compute the SVD of a small matrix $B$: $B =  U \Sigma H^T$.
   }
   \STATE{6. Form the $n \times l$ matrix $V = Q H$, then $K \approx U \Sigma
   V^T$.
   }
\end{algorithmic}
\end{algorithm}

The randomized SVD is much cheaper than the classical SVD. In fact,
the flops count of the classical SVD for matrix $K$
is about $4mn^2 + 8n^3$ \cite{GolubLoan_book13}, while the cost of
Algorithm \ref{algo:RSVD*} is only about $4mnl$
\cite{XiangZou_InverseProb13}.
For the cases where singular values decay rapidly, we can choose $l
\ll m$. The ratio of the costs between RSVD and the classical
SVD is of the order $O(l/n )$  according to the flops.
%
%One improvement is to compute a partial SVD based on Lanczos
%bi-diagonalization \cite{VanHuffelVandewalle_SIAM1991}. This
%procedure will destroy the sparsity or structure of the matrix. Also
%it needs to access the coefficient matrix many times and use the
%BLAS-2 operations, i.e., the matrix-vector multiplications. Hence
%the partial SVD still requires essential CPU times for large-scale
%problems.
%Algorithm \ref{algo:RSVD*} (RSVD) can be essentially less
%expensive than the classical SVD and partial SVD.

We can see that Algorithm \ref{algo:RSVD*} generates an approximate
decomposition $K \approx K Q Q^T = U\Sigma V^T$, where the columns
of $Q$ span approximately the range of $K^T$, or the right singular
vectors.
The RSVD in Algorithm \ref{algo:RSVD*}
%may be viewed as a variant of the randomized SVD in \cite{HalkoMartinssonTropp_SIREV11}, and
was formulated in \cite{XiangZou_InverseProb13}, and can be directly applied for
the matrix $K=P^T A L^\dag Z$ or $K=A L^\# Z$ in the standard system \eqref{Case6ThmSol2} or
\eqref{Case6ThmOnesubprob2} transformed from the system
\eqref{ProbGeneralForm} of general form.
The operations with $K$ involve now the operations with $L^{-1}$,
$L^\dag$, or $L^\#$, which can be implemented efficiently in
many applications.
For example, when $L$ is the discrete Laplacian the actions of inverses
can be done by the algebraic multigrid method efficiently.
For the cases where $L$ is of special structures, such as Hankel or Toeplitz,
there exist many fast solvers for
implementing the operations with $L^{-1}$ or $L^\dag$ \cite{chan07}.
%
%It can be regarded as
%applying the randomized SVD in \cite{HalkoMartinssonTropp_SIREV11}
%on the matrix $K^T$ and achieving $ Q Q^T K^T = V \Sigma U^T$. The
%reason why we use this algorithm in this paper is that for many
%inverse problems the number of rows of the matrix $K$ is much
%smaller the the number of columns.

Suppose that we have an SVD approximation $K \approx U \Sigma V^T$
(by Algorithm \ref{algo:RSVD*}), where $\Sigma$ is diagonal with the
form $\Sigma = \text{diag} ( \sigma_1, \cdots, \sigma_l )$, $U=(u_1,
\cdots, u_l)$ and $V= (v_1, \cdots, v_l)$ are orthonormal matrices.
Then the approximate Tikhonov regularized solution of
\eqref{ProbStandardForm} can be expressed as
$$ x_\mu = \sum_{i=1}^l \frac{\sigma_i^2}{\sigma_i^2+\mu^2} \frac{u_i^T b}{\sigma_i} v_i.$$
The regularization parameter $\mu$ can be determined by several existing
popular methods, such as L-curve, GCV function, or some discrepancy
principles.
If we discard the small diagonal elements in $\Sigma$, we obtain
the truncated SVD (TSVD) of $K$. With an abuse of
notations, we denote the approximate TSVD of $K$ by $K \approx U
\Sigma V^T$, where $\Sigma= \text{diag} (\sigma_1, \cdots,
\sigma_k)$, either of $U$ and $V$ has $k$ orthonormal columns. Then the
approximate TSVD regularized solution $x_k$ is given by
$$ x_k = \sum_{i=1}^k \frac{u_i^T b}{\sigma_i} v_i.$$

\section{Inverse problems of general form and solutions by random generalized SVD}

As we have discussed in the last section, the problem
\eqref{ProbGeneralForm} of general form can be transformed into the
problem of standard form \eqref{ProbStandardForm}, then the
classical or randomized SVD method is applied to seek the
regularized solution. The classical SVD is usually very expensive,
while the randomized SVD method is much cheaper. In this
section we shall discuss an alternative strategy for solving the
problem \eqref{ProbGeneralForm} of general form by using the
generalized SVD (GSVD) of the matrix pair $(A, L)$. But again the
classical GSVD are expensive, so we try to reduce the problem size
and then seek an approximate solution. We will show that the
approximate regularized solution can be achieved by some randomized
algorithms.

\subsection{Regularized solution with exact GSVD}

We consider the problem \eqref{ProbGeneralForm} of general form and
the matrix pair $(A, L)$ with $A \in \mathbb{R}^{m \times n}$, $L \in \mathbb{R}^{p \times
n}$. We assume that $\mathcal{N}(A) \cap \mathcal{N}(L) = \{0\}$, and
$p \geq n$.
The classical generalized SVD (CGSVD) is obtained as follows. We first
perform a QR factorization for the pair $(A, L)$:
$$\begin{bmatrix} A \\ L \end{bmatrix}=
\begin{bmatrix} Q_A
\\ Q_L \end{bmatrix} R, $$
where the matrix $[Q_A; Q_L]$ is column orthonormal. Then the CS
decomposition \cite{GolubLoan_book13,VanLoan_NM85} is applied to
this column orthogonal matrix
$$ \begin{bmatrix} Q_A
\\ Q_L \end{bmatrix} = \begin{bmatrix} U & 0 \\
0 & V \end{bmatrix}  \begin{bmatrix} C \\ S \end{bmatrix} W^T , $$
where $C= \text{diag} (c_1, \cdots, c_n)$, $S= \text{diag}(s_1,
\cdots, s_n)$, $0 \leq c_1 \leq \cdots \leq c_n$, $s_1 \geq \cdots
\geq s_n$, $c_i^2 + s_i^2 = 1$, and $U, V, W$ are orthogonal
matrices with compatible dimensions.
Let $G^{-1} = W^T R$, then we have the classical generalized SVD
(CGSVD) of the matrix pair $(A,L)$ as follows
\cite{GolubLoan_book13}:
\begin{equation}\label{eq:CGSVD}
 A = U C G^{-1}, \qquad L = V S G^{-1}.
\end{equation}
 Let $G=(g_1, \cdots,
g_n)$, $U=(u_1, \cdots, u_n)$ and $V=(v_1, \cdots, v_n)$. Using the
right singular vectors $g_i$, and the two sets of left singular
vectors $u_i$ and $v_i$, we can rewrite the CGSVD for $p \ge n$ as
$$ A g_i = c_i u_i, \quad L g_i = s_i v_i \q \mbox{for} ~~i=1,2, \cdots, n. $$
Now using the above CGSVD, we can find the solution of \eqref{ProbGeneralForm}:
\begin{eqnarray}
x_\mu
%&=& \arg \min_x \left\| \begin{bmatrix} A \\ \mu L \end{bmatrix}
%x -
%\begin{bmatrix} b \\ 0 \end{bmatrix} \right\|
&=& (A^T A + \mu^2 L^TL)^{-1} A^T b
= G (C^T C + \mu^2 S^T S)^{-1} C^T U^T b \nb \\
&=& \sum_{i=1}^{n} \frac{c_i^2}{c_i^2 + \mu^2 s_i^2} \frac{u_i^T
b}{c_i} g_i\,. \l{eq:originalSolutionbyCGSVD}
\end{eqnarray}
%
%The regularized solution can be expressed as
% \begin{equation*} x_\mu = \sum_{i=1}^{n} \frac{c_i^2}{c_i^2 + \mu^2 s_i^2} \frac{u_i^T b}{c_i} g_i
% \end{equation*}

For the case of $p<n$, the generalized singular vectors satisfy the relations:
%\begin{eqnarray*}
%A g_i = c_i u_i,  && L g_i = s_i v_i \q \mbox{for} ~~i= 1, \cdots, p; \\
%A g_i = u_i, && L g_i = 0 \q \mbox{for} ~~i = p+1, \cdots, n.
%\end{eqnarray*}
%
\begin{align*}
&& && &A g_i = c_i u_i,  & &L g_i = s_i v_i  && \mbox{for} ~~i= 1, \cdots, p ; && &&\\
&& && &A g_i = u_i, & &L g_i = 0   && \mbox{for} ~~i = p+1, \cdots,
n. && &&
\end{align*}
Then we can express the regularized solution for $p<n$ as
$$ x_\mu = \sum_{i=1}^p \frac{c_i^2}{c_i^2 + \mu^2 s_i^2} \frac{u_i^T b}{c_i}
g_i + \sum_{i=p+1}^n (u_i^T b) g_i.  $$
%
%For historical reasons, the values $c_1, \cdots, c_p$ are ordered in
%nondecreasing order.
The smoothest GSVD vectors $u_i$, $v_i$ and
$g_i$ are those with $i \approx p$. If $p < n$, then $g_i ~ (i=p+1,
\cdots, n)$ are null vectors of $L$ and therefore very smooth.
We can see that these components
$\{g_i\}_{i=p+1}^n$ are incorporated into the solution directly without any
regularizaiton.

Similarly to TSVD, a truncated version of GSVD can be naturally
extended for the problem \eqref{ProbGeneralForm} of general form.
The truncated GSVD (TGSVD) solution reads
$$ x_k = \sum_{i=p-k+1}^p  \frac{u_i^T b}{c_i} g_i + \sum_{i=p+1}^n (u_i^T b) g_i.  $$

We can use the GSVD based Tikhonov regularization to solve
\eqref{ProbGeneralForm}, but the classical GSVD (CGSVD) above is expensive
and impractical for large-scale problems.
We shall derive a randomized GSVD algorithm that helps us
reduce the large-scale problem size essentially and seek an
approximate regularized solution.
%We will derive this algorithm step by step.

\subsection{Problem size reduction and approximate solution}
\label{sec:size}

Suppose that the matrix $A$ in \eqref{ProbGeneralForm} has the SVD
($m \leq n$) such that $A= U \Sigma V^T$, where $U \in \mathbb{R}^{m
\times m}$ and $V \in \mathbb{R}^{n \times n}$ are unitary matrices,
and $\Sigma \in \mathbb{R}^{m \times n}$ is a diagonal matrix with
the diagonal elements $\sigma_1 \geq \cdots \geq \sigma_m \geq 0 $.
We divide the matrix $\Sigma$ into two parts: $\Sigma = \text{diag}(
\Sigma_1, \Sigma_2 )$, where $\Sigma_1 \in \mathbb{R}^{r \times r}$
and $\Sigma_2 \in \mathbb{R}^{ (m-r) \times (n-r) }$.
Correspondingly we partition the matrices $U$ and $V$ as $U =[U_1,
U_2]$, $V=[ V_1, V_2 ]$, where $U_1 \in \mathbb{R}^{m \times r}$,
$U_2  \in \mathbb{R}^{m \times (m-r)}$, $V_1  \in \mathbb{R}^{n
\times r}$, and $V_2  \in \mathbb{R}^{n \times (n-r)}$.
Then we can split matrix $A$ into two parts:
\begin{equation}\label{eq:svdA}
A = \begin{pmatrix} U_1, U_2 \end{pmatrix} \begin{bmatrix} \Sigma_1
&
\\ & \Sigma_2
\end{bmatrix}
\begin{pmatrix} V_1^T\\ V_2^T \end{pmatrix} =  U_1 \Sigma_1 V_1^T +
U_2 \Sigma_2 V_2^T.
\end{equation}

Suppose that there is a gap among the singular values. The diagonals
of $\Sigma_2$ correspond to the smaller singular values, while the
diagonals of $\Sigma_1$ include the larger ones. Then the matrix can
be approximated by $A \approx U_1 \Sigma_1 V_1^T$. Since the
singular vectors associated with smaller singular values have more
sign changes in their components, we may seek the solution of the
form $x_\mu=V_1 \bar{x}$ to the system \eqref{ProbGeneralForm}
and come to solve the following problem of the reduced size:
\begin{equation}\label{ProblemReduced_GSVD}
\min_{\bar{x}} ||A V_1 \bar{x} - b||^2 + \mu^2 ||L V_1 \bar{x}||^2 .
\end{equation}
It is equivalent to
$$
\min_{\bar{x}} \left\|  \begin{bmatrix} A  \\
\mu L
\end{bmatrix} V_1 \bar{x} - \begin{bmatrix} b \\ 0 \end{bmatrix}
\right\| .
$$
 Since $\mathcal{N}(A) \cap \mathcal{N}(L)
= \{0\}$, $\begin{bmatrix} A  \\ \mu L  \end{bmatrix}$  is of
full column rank. As $V_1$ is column orthogonal, so $\begin{bmatrix} A  \\
\mu L \end{bmatrix} V_1 $ is also of full column rank. Hence the
reduced problem \eqref{ProblemReduced_GSVD} has a unique solution:

\begin{equation}
\l{eq:approximate_soln}
\bar{x}_\textsc{ls} =  \left(V_1^T (A^T A + \mu^2 L^T L) V_1
\right)^{-1} V_1^T A^T b.
\end{equation}

Note that for this approximate regularized solution to \eqref{ProblemReduced_GSVD} we only need to
work with the matrix pair $( A V_1, L V_1 )$ with the size of $m
\times r$ and $p \times r$ respectively, while the original matrix
pair $(A, L)$ is of size $m \times n$ and $p \times n$ respectively.
We often take $r\ll n$, so the size of the approximate system \eqref{ProblemReduced_GSVD}
is essentially smaller than the original one \eqref{ProbGeneralForm}.
We shall only work on the reduced problem, hence the memory requirement
and CPU time can be significantly reduced.

Next, we shall compare the approximate solution \eqref{eq:approximate_soln} to the reduced system
\eqref{ProblemReduced_GSVD}
with the exact solution to the system \eqref{ProbGeneralForm} of general form.
To do so, we represent
$\bar{x}_\textsc{ls}$ in terms of the SVD \eqref{eq:svdA} of $A$.
Define $\bar{\Sigma}_2^2 = \Sigma_2^T \Sigma_2 = \text{diag}
(\sigma_{r+1}^2, \cdots, \sigma_m^2, 0, \cdots, 0) \in
\mathbb{R}^{(n-r) \times (n-r)}$.
Now a direct computing yields that
\begin{eqnarray} A^T A + \mu^2 L^T L & = & V \left( \begin{bmatrix} \Sigma_1^2 & \nb\\
& \bar{\Sigma}_2^2 \end{bmatrix} + \mu^2 (L V)^T (LV) \right) V^T
\\ & = & V \begin{bmatrix} F & B \\ B^T & D \end{bmatrix} V^T,
\l{eq:schur}
\end{eqnarray}
 where $F$, $B$ and $D$ are given by
\begin{eqnarray*} F &=& \Sigma_1^2 + \mu^2 (L V_1)^T  L V_1  \in \mathbb{R}^{r \times r}, \\
B &=& \mu^2 (L V_1)^T L V_2 \in \mathbb{R}^{r \times (n-r)}, \\
D &=& \bar{\Sigma}_2^2 + \mu^2 (L V_2)^T  L V_2 \in
\mathbb{R}^{(n-r)\times (n-r)}.
\end{eqnarray*}
It is easy to see that $F$ is nonsingular, and we can write the
solution of \eqref{ProblemReduced_GSVD} as follows:
\begin{eqnarray*}
\bar{x}_\textsc{ls} %&=& \left(V_1^T (A^T A + \mu^2 L^T L) V_1 \right)^{-1} V_1^T A^T b   \\
&=& \left(V_1^T V \begin{bmatrix} F & B \\ B^T & D \end{bmatrix} V^T
V_1 \right)^{-1} V_1^T V \begin{bmatrix} \Sigma_1 &  \\   &
\Sigma_2^T \end{bmatrix} U^T b \\
&=& \left( [I,0] \begin{bmatrix} F & B \\ B^T & D \end{bmatrix}
\begin{bmatrix} I \\ 0 \end{bmatrix}
 \right)^{-1} [I,0] \begin{bmatrix} \Sigma_1 &  \\   & \Sigma_2^T
\end{bmatrix} U^T b \\
&=& F^{-1} \Sigma_1 U_1^T b .
\end{eqnarray*}
Then we obtain an approximate regularized solution of \eqref{ProbGeneralForm}:
\begin{equation}\label{ProblemReduced_GSVD_ApproxSol}
 x_\mu \approx V_1
\bar{x}_\textsc{ls} = V_1 F^{-1} \Sigma_1 U_1^T b.
\end{equation}

For the special case where $L=I$, the formula
\eqref{ProblemReduced_GSVD_ApproxSol} is simplified as
$$ x_\mu \approx \sum_{i=1}^r \frac{\sigma_i^2}{\sigma_i^2+\mu^2} \frac{u_i^T b}{\sigma_i} v_i,$$
which is further reduced for $\mu=0$ as
$$ x_\mu \approx \sum_{i=1}^r \frac{u_i^T b}{\sigma_i} v_i.$$
This can be regarded as a regularized solution by truncation.

Next, we analyse the difference between the approximate solution
\eqref{ProblemReduced_GSVD_ApproxSol} and the original exact
solution \eqref{eq:originalSolutionbyCGSVD}.
%\footnote{This
%paragraph is quite important. But it is not straightforward to
%understand which difference you are discussing, please make this
%paragraph more clear and point out the difference explicitly.
%\textcolor{blue}{Sentences are added before this place. Nearly 16
%lines are added.} }
This process will show us how to obtain this approximate
solution alternatively.
For this purpose, we define the Schur complement $S=D-B^T F^{-1} B$
associated with the governing matrix in \eqref{eq:schur}. It is easy
to check that
$$ y^T S y = y^T  \bar{\Sigma}_2^2 y + \mu^2 (L V_2 y)^T [I-
\mu^2 L V_1 F^{-1} (L V_1)^T] L V_2 y. $$
%$L V_2 y \neq 0$???.
Since $\mu$ is a small parameter, it is reasonable to assume that $
y^T S y > 0$ for any $y \neq 0$, hence $S$ is nonsingular.
%After some approximations, the term involving $S$ will disappear  in our
%final expression, but we need it temporarily  in the following
%derivation.
%
We can verify that
$$ \begin{bmatrix} F & B \\ B^T & D \end{bmatrix}^{-1}  = \begin{bmatrix} M_{11} & M_{12} \\ M_{21} &
S^{-1}
\end{bmatrix},$$
where
$%\begin{eqnarray*}
M_{11} = F^{-1} + F^{-1} B S^{-1} B^T F^{-1}, $ and $ M_{21} =
M_{12}^T = - S^{-1} B^T F^{-1}
$%\end{eqnarray*}
.
We can write the solution of \eqref{ProbGeneralForm} (see
\eqref{eq:originalSolutionbyCGSVD}) as
\begin{eqnarray}
x_\mu & = & ( A^T A + \mu^2 L^T L )^{-1} A^T b \nb\\
 & = & [V_1, V_2] \begin{bmatrix} M_{11} & M_{12} \\ M_{21} & S^{-1}
\end{bmatrix} \begin{bmatrix} \Sigma_1 &  \nb\\   & \Sigma_2^T
\end{bmatrix} \begin{bmatrix} U_1^T \\ U_2^T \end{bmatrix} b \nb\\
& = & V_1 M_{11} \Sigma_1 U_1^T b +  V_2 M_{21} \Sigma_1 U_1^T b  +
(V_1 M_{12} + V_2 S^{-1} ) \Sigma_2^T U_2^T b  \nb \\
& = & V_1 (I + F^{-1} B S^{-1} B^T ) F^{-1} \Sigma_1 U_1^T b
 - V_2 S^{-1} B^T F^{-1} \Sigma_1 U_1^T b   \nb \\ & &
 + ( V_2 - V_1
F^{-1} B ) S^{-1} \Sigma_2^T U_2^T b . \l{eq:sigma2}
\end{eqnarray}
%\footnote{I am considering how to reuse the 2nd term: $(V_2 V_2^T) (
%V_2 S^{-1} B^T V_1^T ) (V_1^T F^{-1} \Sigma_1 U_1^T b) $. I am
%feeling that it would be useful ...}

We shall show how to achieve the approximate regularized solution
\eqref{ProblemReduced_GSVD_ApproxSol} from the exact solution
\eqref{eq:sigma2} by truncation.
Suppose that there exists a gap between $\Sigma_1$ and $\Sigma_2$.
And we further assume that $\Sigma_1=O(1)$, and
$\Sigma_2=O(\epsilon)$ or $o(\epsilon)$. Then the third term in
\eqref{eq:sigma2}  may be ignored, leading to
$$x_\mu \approx V_1 (I + F^{-1} B S^{-1} B^T ) F^{-1} \Sigma_1 U_1^T b
 - V_2 S^{-1} B^T F^{-1} \Sigma_1 U_1^T b   . $$
 Note that $V_2$ is associated with the
smallest singular values and its columns are highly oscillatory and
have more frequent sign changes. If we drop the second term involving $V_2$
in \eqref{eq:sigma2} like we do with TSVD, we then derive
$$x_\mu \approx V_1 (I + F^{-1} B S^{-1} B^T ) F^{-1} \Sigma_1 U_1^T b. $$
The expected regularization parameters are usually small, and it is
reasonable to assume that $\mu=O(\epsilon)$. Then we can see that
$B=O(\epsilon^2)$, $F^{-1}=O(1)$, $S=O(\epsilon^2)$, and $F^{-1} B
S^{-1} B^T = O(\epsilon^2)$. If we further  omit the %$O (\epsilon)$
higher order terms in the above expression, then
we obtain the following approximate solution
\begin{equation}\l{eq:approximate2}
x_\mu \approx V_1 F^{-1} \Sigma_1 U_1^T b ,
\end{equation}
which is exactly the approximate solution $V_1 \bar{x}_\textsc{ls}$
in \eqref{ProblemReduced_GSVD_ApproxSol} that we will use in the
following randomized algorithm.

We end this subsection with some remarks on a few important cases where we may
further simplify the representation of solutions.

{\bf Case 1: $L=I$}. For this simplest case, we can check that
$B=0$, $F=\Sigma_1^2 + \mu^2 I$, and $S = D=\bar{\Sigma}_2^2 + \mu^2
I$. Then the exact solution \eqref{eq:sigma2} of \eqref{ProbGeneralForm}  can be expressed
by
$$ x_\mu  = V_1   F^{-1} \Sigma_1 U_1^T b    +   V_2    S^{-1} \Sigma_2^T U_2^T b , $$
which is the same as $(A^T A + \mu^2 I)^{-1} A^T b$ by using
the SVD \eqref{eq:svdA} of $A$. If we chop off the second oscillation term, we have $ x_\mu
\approx V_1 F^{-1} \Sigma_1 U_1^T b $ in \eqref{eq:approximate2}.

{\bf Case 2: span$\{ V_1 \} = \mathcal{R} (A^T)$}. In this case we
know $\Sigma_2$ is simply a zero matrix, and the last term in
\eqref{eq:sigma2} vanishes. Then the exact solution \eqref{eq:sigma2}  of
\eqref{ProbGeneralForm} reduces to
$$
x_\mu = V_1 M_{11} \Sigma_1 U_1^T b + V_2 M_{21} \Sigma_1 U_1^T b.$$
Note that the second term above lies in $\mathcal{N} (A)$ and is
perpendicular to the first term. Hence the least-squares solution
with minimum norm is given by
$$ x_\mu = V_1 M_{11} \Sigma_1 U_1^T b %= V_1 (I + F^{-1} B S^{-1} B^T ) F^{-1} \Sigma_1 U_1^T b
.$$ As we have seen before, this term can be further approximated as
$ x_\mu  \approx V_1 F^{-1} \Sigma_1 U_1^T b $.

{\bf Case 3: span$\{ V_2 \} \subseteq \mathcal{N} (L) \subseteq
 \mathcal{N} (A) $}. For this case, we can see that all the matrices
 $B$, $\Sigma_2$, $S$ and $D$ vanish. Then the exact solution \eqref{eq:sigma2} of \eqref{ProbGeneralForm}
 is simplified to
 $$ x_\mu =  V_1   F^{-1} \Sigma_1 U_1^T b,$$
which is exactly the same as $(A^T A + \mu^2 L^TL)^\dag A^T b$.

%Case 3'. L=0

%-------------------------

\subsection{Regularized solution by the randomized GSVD}

As discussed in section\,\ref{sec:size}, we need to obtain a good approximation
$V_1$ required in the system \eqref{ProblemReduced_GSVD}, in order to
reduce the size of the problem \eqref{ProbGeneralForm}
of general form and find its approximate solution. If we directly perform the SVD
of $A$ and choose $V_1$ from its right singular vectors, it will be very expensive
and impractical for large-scale systems.
We now seek an economic way to obtain a good approximation $V_1$.
%required in the system \eqref{ProblemReduced_GSVD} of reduced size.
 Suppose that we have an approximate SVD of $A$, that is, $A
\approx U \Sigma V^T$, where $U \in \mathbb{R}^{m \times l}$ and $V
\in \mathbb{R}^{n \times l}$ are orthonormal. There is an
abuse of notation $U$ and $V$ here in GSVD of $(A,L)$, but they can be
differentiated from the context. The approximate SVD of $A$ can be achieved by RSVD, i.e.,
Algorithm \ref{algo:RSVD*} with $K$ replaced by $A$. We write this approximation as $A \approx U
\Sigma \tilde V_1^T$, then we can seek the solution of the form $x =\tilde V_1 \widehat{x}
$, where the ${n \times l}$ column orthogonal matrix  $\tilde V_1$ forms the
approximate right singular vectors of $A$. Now the transformed problem \eqref{ProblemReduced_GSVD}
reads as follows:
\begin{equation}\l{eq:reduced}
\min_{\widehat{x}}  ||A \tilde V_1 \widehat{x} - b||^2 + \mu^2 ||L \tilde V_1 \widehat{x}||^2 .
\end{equation}
Define $\widehat{A}= A \tilde V_1$ and $\widehat{L}= L \tilde V_1$.  This problem
still has the general form

\begin{equation}\label{ProblemReducedbyRGSVD}
\min_{\widehat{x}}  || \widehat{A} \widehat{x} - b||^2 + \mu^2 ||
\widehat{L} \widehat{x}||^2 ,
\end{equation}
where $\widehat{A} \in \mathbb{R}^{m \times l}$ and $\widehat{L} \in
\mathbb{R}^{p \times l}$.
Since $\mathcal{N}(A) \cap \mathcal{N}(L) = \{0\}$ and $\tilde V_1$ is column orthogonal,
then
$[\widehat{A}^T,  \mu \widehat{L}^T]^T$ is  of full column rank. Hence the
reduced problem \eqref{ProblemReducedbyRGSVD} has a unique
solution.
 But the matrix pair $(\widehat{A},
\widehat{L})$ is of much smaller size compared with the original matrix
pair $( A, L )$, and we can easily apply the classical GSVD (CGSVD)
to this matrix pair. The solution procedure is summarized in
Algorithm \ref{algo:RGSVDregu}.

\begin{algorithm}
\caption{ Approximate regularized solution by RGSVD.   }
\label{algo:RGSVDregu}
\begin{algorithmic}
   \STATE{1. Use RSVD, i.e., Algorithm \ref{algo:RSVD*},  to obtain
   }
   \STATE{\qquad $A \approx U \Sigma \tilde V_1^T$  with $U \in \mathbb{R}^{m\times l}$, $\Sigma \in \mathbb{R}^{l \times l}$ and $\tilde V_1 \in \mathbb{R}^{n \times l}$.
   }
   \STATE{2. Apply CGSVD to the matrix stencil $(A \tilde V_1, L \tilde V_1)$ of smaller size. % by calling QR and CS decompositions.
   }
   \STATE{3. Solve the Tikhonov regularization problem \eqref{eq:reduced}.
   }
   %\STATE{\qquad $\min_{\widehat{x}} \{ || A \tilde V_1 \widehat{x} - b||^2 + \mu^2 || L \tilde V_1 \widehat{x}||^2 \} $.
   \STATE{4. Set $x = \tilde V_1 \widehat{x}$.
   }
\end{algorithmic}
\end{algorithm}

%-----------------------------------------------

Suppose that the CGSVD of the matrix pair $(\widehat{A},
\widehat{L})$ has the similar form as \eqref{eq:CGSVD},
$$ \widehat{A} = \widehat{U} \widehat{C} \widehat{G}^{-1}, \qquad \widehat{L} = \widehat{V} \widehat{S} \widehat{G}^{-1}. $$
Then the solution $\widehat{x}$ of \eqref{ProblemReducedbyRGSVD} can
be expressed by
$$\widehat{x}_\textsc{ls} = \widehat{G}  (\widehat{C}^T \widehat{C} + \mu^2 \widehat{S}^T
\widehat{S})^{-1} \widehat{C}^T \widehat{U}^T b .$$
And this gives us an approximate solution of the original problem \eqref{ProbGeneralForm}:
\bb \l{eq:new_approximate_soln}
x_\mu \approx \tilde V_1 \widehat{x}_\textsc{ls} = \tilde V_1 \widehat{G}  (\widehat{C}^T \widehat{C} + \mu^2 \widehat{S}^T
\widehat{S})^{-1} \widehat{C}^T \widehat{U}^T b .
\ee

In the following we will demonstrate that this approximation can be regarded
as the least-squares solution with minimum norm of a nearby problem.
To this aim, we define
\begin{equation}\label{Decompositon_RGSVD}
 \begin{bmatrix} \widetilde{A} \\ \widetilde{L} \end{bmatrix} := \begin{bmatrix} A \\ L \end{bmatrix} \tilde V_1 \tilde V_1^T = \begin{bmatrix} \widehat{A} \\ \widehat{L} \end{bmatrix} \tilde V_1^T
= \begin{bmatrix} \widehat{U} \widehat{C} \widehat{G}^{-1} \tilde V_1^T \\
\widehat{V} \widehat{S} \widehat{G}^{-1} \tilde V_1^T \end{bmatrix}.
\end{equation}
%This can be regarded as an approximation decomposition of
%the matrix pair $(A, L)$, and we call this the randomized version of
%generalized singular value decomposition (RGSVD).

%
\begin{lemma}  The approximate solution \eqref{eq:new_approximate_soln} is the solution to the following problem
with minimum norm:
\begin{equation}\label{ProblemNearby_afterRGSVD}
\min_x || \widetilde{A}  x - b||^2 + \mu^2 || \widetilde{L} x ||^2\,.
\end{equation}
\end{lemma}
%where $\widetilde{A}=A\tilde V_1\tilde V_1^T $ and $\widetilde{L}=L\tilde V_1\tilde V_1^T $.

\no {\bf Proof}.
%Under the assumption $\mathcal{N}(A) \cap \mathcal{N}(L) = \{0\}$,
%the matrix $[A;L]$ is of full column rank, and the original problem
%\eqref{ProbGeneralForm} has a unique solution. This is not the case
%for the matrix $[\widetilde{A}; \widetilde{L}]$ any more.
We can
check that $\mathcal{N}(\widetilde{L}) = \mathcal{N}(\widetilde{A})=
span\{\tilde V_1\}^\perp$, which is included in Case 1
in section~\ref{sec:special_cases}, from where we know minimisation \eqref{ProblemNearby_afterRGSVD}
has the least-squares solution with
minimum norm:
$$
\widetilde{x}_\mu = ( \widetilde{A}^T \widetilde{A} + \mu^2 \widetilde{L}^T \widetilde{L} )^\dag \widetilde{A}^T b.
$$
Using the decomposition \eqref{Decompositon_RGSVD} and the property
of Moore-Penrose inverse, this regularized solution can be further expressed as
\begin{eqnarray*}
\widetilde{x}_\mu %& := &  ( \widetilde{A}^T \widetilde{A} + \mu^2 \widetilde{L}^T \widetilde{L} )^\dag \widetilde{A}^T b \nb \\
&=& [\tilde V_1 \widehat{G}^{-T} (\widehat{C}^T \widehat{C} + \mu^2 \widehat{S}^T \widehat{S})  \widehat{G}^{-1} \tilde V_1^T]^\dag \widetilde{A}^T b \nb \\
&=& \tilde V_1 \widehat{G}  (\widehat{C}^T \widehat{C} + \mu^2 \widehat{S}^T \widehat{S})^{-1} \widehat{G}^T \tilde V_1^T  \widetilde{A}^T b \nb \\
&=& \tilde V_1 \widehat{G}  (\widehat{C}^T \widehat{C} + \mu^2 \widehat{S}^T
\widehat{S})^{-1} \widehat{C}^T \widehat{U}^T b \,.
\l{eq:approxSolutionbyRGSVD}
\end{eqnarray*}
This is exactly the solution \eqref{eq:new_approximate_soln}.
\qed

%Obviously it has the similarity with the expression
%\eqref{eq:originalSolutionbyCGSVD}, and it is the approximate
%regularized solution based on RGSVD.

In summary, for the least-squares solution of general form
\eqref{ProbGeneralForm}, we can approximate it by solving the problem
\eqref{ProblemReducedbyRGSVD} and set $x=\tilde V_1 \widehat{x}$, which is
equivalent to the least-squares solution of
\eqref{ProblemNearby_afterRGSVD} with minimum norm. Note that the
size of the matrix pair $(\widehat{A},\widehat{L})$ is much smaller
than the matrix pairs $(A,L)$, and we need to work only on the matrix pair $(\widehat{A},\widehat{L})$
in practice.
Using Algorithm \ref{algo:RGSVDregu}, we can obtain an approximate
solution \eqref{ProblemReduced_GSVD_ApproxSol} based on the GSVD on
the matrix pair $(\widehat{A},\widehat{L})$ with $\widehat{A} \in
\mathbb{R}^{m \times l}$ and $\widehat{L} \in \mathbb{R}^{p \times
l}$, while the original problem of general form has the regularized
solution \eqref{eq:originalSolutionbyCGSVD} when CGSVD is applied
directly on the matrix pair $(A,L)$.
Furthermore, we know the approximate solution is spanned by the
columns of $\tilde V_1$. So if we have some a priori information about the
solution, we may incorporate it into $\tilde V_1$ by orthogonalization with
its columns.

%TGSVD, page 51.

Similarly to the truncated SVD (TSVD), we can use the truncated GSVD
(TGSVD) to seek the regularized solution \cite{Hansen_book98}. That
is, after using Algorithm \ref{algo:RGSVDregu}, we chop off those
smallest singular values to achieve a truncated version of RGSVD.
Then we use the L-curve, GCV, or other discrepancy principles to
determine the regularization parameter. This procedure is similar to
Tikhonov regularization, and we do not go into further details about
this.

\section{Numerical experiments}

For the system \eqref{ProbStandardForm} of standard form, we
studied RSVD and tested various linear inverse problems in
\cite{XiangZou_InverseProb13}. In this section, we shall focus on
the system \eqref{ProbGeneralForm} of general form directly and test the newly proposed
Algorithm \ref{algo:RGSVDregu} with examples from different linear
inverse problems to illustrate the performance of the algorithm. As
we shall see, we can essentially reduce the problem size by using
RGSVD, but still obtain reasonably good approximate regularized
solutions. The following tests are done using \textsc{Matlab} R2012a
in a laptop with Intel(R) Core(TM) i5 CPU M480 @2.67G.

In this section we will test some examples from Regularization Tools
\cite{Hansen_regutool}. Most of the cases are related to the
Fredholm integral equation of the first kind,
\begin{equation} \label{Fredholm1st}
\int_a^b K(s,t) x(t) dt = g(s), \qquad c \leq s \leq d,
\end{equation}
where $K$ is the square integrable kernel function. In the tests the
kernel $K$ and the solution $x$ are given and discretized to yield
the matrix $A$ and the vector $x$, then the discrete right-hand side
is determined by $b=Ax$. Matrix $L$ in the general form
\eqref{ProbGeneralForm} is the discrete approximation to  the first
or second order differential operator. In particular we will try
the following choices:
$$ L = \begin{pmatrix} 1 & -1 & & & \\   & 1 & -1 & &  \\   &  & \ddots & \ddots & \\   &   &   & 1 & -1 \end{pmatrix}  \quad \text{and} \quad L = \begin{pmatrix} 1 & -2 & 1 & & & \\   & 1 & -2 & 1 & & \\   &  & \ddots & \ddots & \ddots &\\   &   &   &1 & -2 & 1 \end{pmatrix} \,. $$
We will write them by $L=tridiag(1,-1)$  and $L=tridiag(1,-2,1)$ respectively.

In all our numerical experiments, the observation data $b^\delta$ is
generated from the exact data $b$ by adding the noise in the form
$$b^\delta = b +  \delta  \frac{||b||}{||s||} s = b +  \varepsilon  \frac{s}{||s||}, $$
where $s$ is a random vector,
$s = {randn}(n,1)$ if not specified otherwise, %NOS, eg. Sauf indication contraire, veuillez s'il vous plait, utiliser ma carte de crédit pour toutes les commandes
$\varepsilon = \delta ||b||$
is the so-called noise level, %\cite{BrezinskiRodriguezSeatzu_NA08},
and $\delta$ is the relative noise level \cite{Hansen_book2010}.
In our numerical tests, we choose the relative noise level
$\delta=$1e-4 as in \cite{Hansen_NLAA13}.

In the following tests, we will compare the computational time and
solution accuracy of the RGSVD based Algorithm \ref{algo:RGSVDregu}
with other traditional methods by using CSVD ($L=I$), or CGSVD. The
sampling size is simply chosen to be $l=50$ in the randomized
algorithms if not specified otherwise.
The total CPU times $T$ (in seconds) for seeking regularized
solutions are recorded, which includes the time for  (approximate)
matrix decomposition, regularization parameter determination, and
Tikhonov regularized solution.
In the subsequent numerical tables we shall use the following notation:
$\mu$  stands for the regularization parameter determined
by GCV function, ${err}$ for the relative error
$||x_{\mu}-x||_2/||x||_2$ of the Tikhonov regularized solution
$x_\mu$ to the exact solution $x$,
$T(s)$ is the total CPU time for seeking regularized solution (in
seconds) and $n$ for the problem size.
%The following test cases come from \cite{Hansen_regutool}.

{\bf Example 1} (\textsc{Shaw}). This is a one-dimensional model of an
image reconstruction problem. It arises from discretization of the integral
equation \eqref{Fredholm1st} with the kernel $K$
being the point spread function for an infinitely long slit:
$$ K(s, t) = \left( \cos s + \cos t \right)^2 \left(  {\sin
u} / {u} \right)^2, \quad u= \pi ( \sin s + \sin t ). $$
The exact solution is given by $x(t)=a_1 \exp( -c_1(t-t_1)^2 ) + a_2
\exp( -c_2(t-t_2)^2 )  $, where the parameters $a_1, a_2$, etc., are
constants chosen to give two different humps \cite{Hansen_regutool}.

\begin{figure}[htbp] \centering
\includegraphics[width=\textwidth]{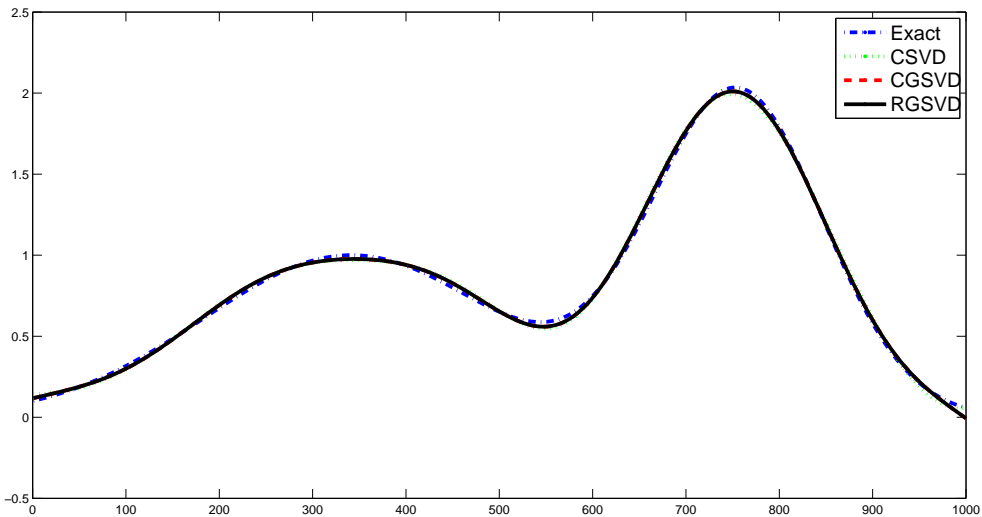}
\caption{  \textsc{Shaw} matrix of size $n=1000$. }
\label{subFig_Shaw1000}
\end{figure}

The computed regularized solutions are shown in Figure
\ref{subFig_Shaw1000}. We can see that the result of the new algorithm RGSVD (black
solid line) is as good as that
of CGSVD (red dashed line). %\footnote{Fig 1(a) and Fig 2(a) were plotted in July 2011.}
But using RGSVD, we only need to work on a much smaller matrix pair.
The total computation times, the regularization parameters and the
relative errors of the computed solutions are given in Table
\ref{tab:Shaw}.
%\footnote{Why do we need to show the results for CSVD ?
%CSVD does not apply to the general form (2), so what problem do you try
%CSVD here ? It is very confusing with CSVD. Please either explain it more clearly
%or remove those results. Also what operator $L$ is used here ?}.
For the problem of size $n=2000$, the new method is
about 100 times faster than the traditional method using CGSVD. When
the problem size is larger, the advantage of the new method is more
obvious. It is clear that the solution of general form with $L=tridiag(1,-2,1)$ %\footnote{What operator $L$ is used here ?Please specify.}
is better than that of the standard form ($L=I$) according to the
relative errors. Moreover, the accuracy of the regularized
solution based on RGSVD is comparable with that of the solution via
CGSVD, but the computation time is much less since we need only to work on a
problem with much smaller size. These observations
also apply to other subsequent testing cases and will not be repeated any more.

\newsavebox{\tableboxShaw}
\begin{lrbox}{\tableboxShaw}

\begin{tabular}{cccccccccccc}\hline
& \multicolumn{3}{c}{CSVD} & & \multicolumn{3}{c}{ CGSVD } & & \multicolumn{3}{c}{ RGSVD }\\
 \cline{2-4} \cline{6-8} \cline{10-12}
 $n$  & $T(s)$ & $\mu$ & $err$ & &   $T(s)$ & $\mu$ & $err$  & &  $T(s)$ & $\mu$ & $err$ \\
\hline

                500   & 0.46 & 4.59E-04 & 3.22E-02 & & 0.47  & 1.00E-01 &  2.16E-02  & & 0.06  & 1.00E-01 &  2.16E-02 \\
                1000  & 1.49 & 2.59E-04 & 2.74E-02 & & 2.57  & 3.34E-01 &  1.89E-02  & & 0.09  & 3.34E-01 &  1.89E-02 \\
                2000  & 9.73 & 2.59E-04 & 3.54E-02 & & 20.6  & 1.21E-00 &  3.02E-02  & & 0.19  & 1.21E-00 &  3.02E-02 \\

\hline
\end{tabular}

\end{lrbox}

\begin{table}
 \centering \scalebox{0.8}{\usebox{\tableboxShaw}}
 \caption{\label{tab:Shaw} Comparison of the CPU time and
solution accuracy among CSVD($L=I$), CGSVD and RGSVD.
%The testing matrix is \textsc{Shaw} with size of $n$.
%The relative noise level is $\delta=$1e-4, and the noise
%levels $\varepsilon$= 0.0052, 0.0074, 0.0104 respectively.
%The
%sampling size is chosen to be $l=50$ for the randomized algorithm.
}
\end{table}

{\bf Example 2} (\textsc{I\_laplace}). This test problem is the
inverse Laplace transformation, a Fredholm first kind integral
equation, discretized by Gauss-Laguerre quadrature. The kernel $K$
is given by $K(s,t) = exp(-st)$.

Regularization Tools \cite{Hansen_regutool} provides the test
problem \textsc{I\_laplace}($n,eg$), with $n$ being the matrix size,
and $eg=1, 2, 3$ or $4$ corresponds to four different examples.
The test problem \textsc{I\_laplace}($n$, 1) has
 the exact solution $x(t) = exp (-t/2)$, and \textsc{I\_laplace}($n$, 3) gives $x(t) =
 t^2 exp(-t/2)$.
%
%We use $L = tridiag(-1,1)$, the discrete approximation to the derivative operator of order 1.
%
For these two cases, the solutions are smooth. The new method
works very well with regularisation $L=tridiag(1,-1)$ and small sample size
$l=50$. The accuracy is comparable or even better than that of the classical
method using CGSVD, but the CPU times and memory requirements are
essentially reduced. We do not report the numerical results
for these two relatively easy cases, but focus on the other two cases, $eg=2$, $4$,
which are more difficult due to the sudden change and strong
discontinuity in the solutions. For these two cases,
the regularization in general form \eqref{ProbGeneralForm} is
necessary to ensure a meaningful numerical solution.

%-----------------------------------------------------------
%-----I_Laplace(n,2)

\begin{figure}[htbp] \centering
\includegraphics[width=\textwidth]{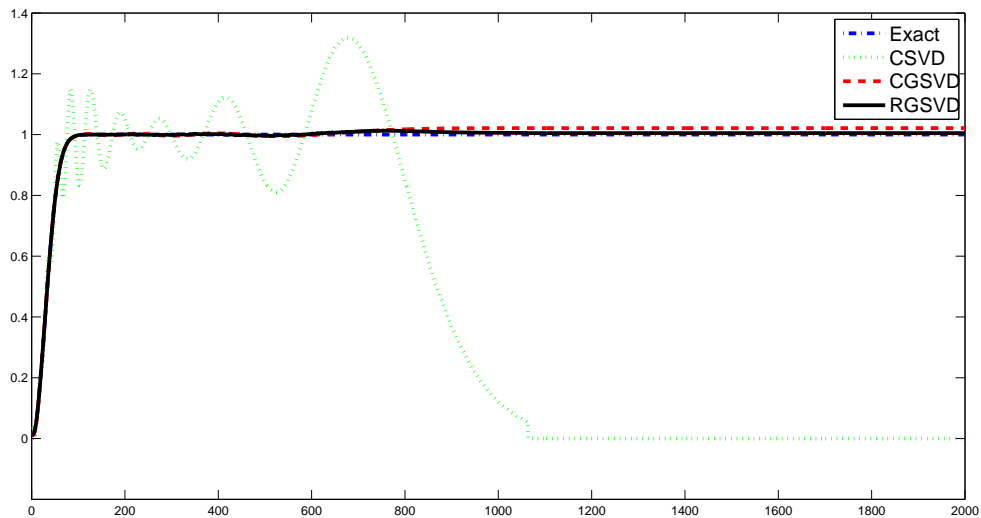}
\caption{  \textsc{I\_laplace}($n$,2) matrix of size $n=2000$. }
\label{fig_ilaplace2_2000}
\end{figure}

First for the test problem \textsc{I\_laplace}($n$, 2), the exact
solution is $x(t)=1-exp(-t/2)$, which has a horizontal asymptote.
For this problem, the identity regularisation $L=I$ can not give a
good reconstruction \cite{Hansen_NLAA13} (see Figure \ref{fig_ilaplace2_2000}), and the regularisation
$L=tridiag(1,-2,1)$ does not work well either. Instead, the
regularisation $L=tridiag(1,-1)$ is very effective to capture the
rapid change in the solution.

Let $e$ be the vector of all ones, i.e., $e=ones(n,1)$. Clearly $e$
is the basis vector of the null space of operator $L=tridiag(1,-1)$.
This suggests us to incorporate a constant mode into the matrix $W$
in Algorithm \ref{algo:RGSVDregu}.
%
%Define the vector $e$ as a vector of all ones, i.e., $e=ones(n,1)$,
%\textcolor{blue}{the basis vector for the null space of
%$L=tridiag(1,-1)$.}
Suppose we have the approximate SVD of the
matrix $A$ by randomized algorithm. That is, $A \approx U \Sigma
W^T$. Let $w  = ( I - W W^T )e$, then we enlarge the matrix $W$ by adding one more column
vector, namely ${w}/{||w||}$.
This is equivalent to finding the orthogonal projection onto $span\{W\}^\perp$,
then adding it to matrix $W$ as a new column.

For most of our cases, the sample size can be as small as 50. But
for this difficult case, we need to use larger size of samples. Even with
larger sample size, the computational time of the new method is
still much less than the traditional method, and the approximate solution is still quite
accurate, actually the accuracy of the new method is much better
than the traditional method; see Figure \ref{fig_ilaplace2_2000} and
Table \ref{tab:ILaplaceTWO} for more details.

\begin{figure}[htbp] \centering
\includegraphics[width=\textwidth]{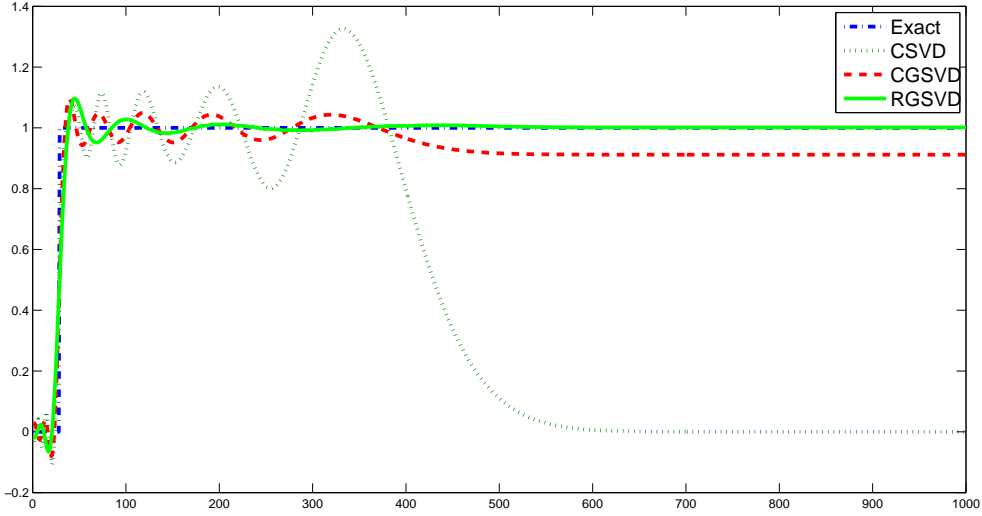}
\caption{  \textsc{I\_laplace}($n$,4) matrix of size $n=1000$. }
\label{fig_ilaplace4_1000}
\end{figure}

We observe from all the examples we have tested, our method working
on a much smaller-sized problem is often more robust and stable than
the classical method using CGSVD. This is also clearly observed in
our test with \textsc{I\_laplace}($n$, 4) from Regularization Tools
\cite{Hansen_regutool}. The exact solution to the problem has a big
jump:
$$ x(t) = \left\{ \begin{matrix} 0, \quad t \leq 2 ; \\ 1, \quad t>2 . \end{matrix} \right. $$
We choose the regularisation
$L=tridiag(1,-1)$.
%In our experiments,
%our method using RGSVD is almost always better than the method using
%CGSVD. That is, our method is more robust, even though we work on a
%smaller size problem (see Fig. \ref{fig_ilaplace4_1000}). In Table
%\ref{tab:ILaplaceFOUR}, the accuracy of our method is always better
%than that of the traditional method. Our method can give somewhat
%accurate solution even when the traditional method fails.
%
%
We have observed that the traditional method using CGSVD fails mostly when we run the same test with
different random noise in the data (right-hand side), but the new method using
RGSVD always succeeds and achieves much better accuracy and requires much less time
than the traditional method, even though we work
on an approximate problem of much smaller size;
see more details in Figure\,\ref{fig_ilaplace4_1000} and Table\,\ref{tab:ILaplaceFOUR}. %\ref{tab:gravity}

We remark that the truncated version of RGSVD also works
quite well for the examples. Since the results are similar to the ones by Tikhonov regularization,
we do not show any results by TGSVD in this work.

\newsavebox{\tableboxILaplaceTWO}
\begin{lrbox}{\tableboxILaplaceTWO}

\begin{tabular}{cccccccccccc}\hline
& \multicolumn{3}{c}{CSVD} & & \multicolumn{3}{c}{ CGSVD } & & \multicolumn{3}{c}{ RGSVD }\\
 \cline{2-4} \cline{6-8} \cline{10-12}
 $n$  & $T(s)$ & $\mu$ & $err$ & &   $T(s)$ & $\mu$ & $err$  & &  $T(s)$ & $\mu$ & $err$ \\
\hline

                500   & 0.44 & 2.15E-04 & 7.59E-01 & & 0.47  & 2.34E-02 &  5.10E-03  & & 0.21  & 1.60E-01 &  1.50E-03 \\
                1000  & 1.56 & 1.18E-04 & 7.56E-01 & & 3.59  & 1.26E-01 &  7.40E-03  & & 1.34  & 2.47E-01 &  1.80E-03 \\
                2000  & 8.76 & 2.43E-04 & 7.45E-01 & & 18.2  & 1.89E-01 &  1.71E-02  & & 5.98  & 4.30E-01 &  5.50E-03 \\

\hline
\end{tabular}

\end{lrbox}

\begin{table}
 \centering \scalebox{0.8}{\usebox{\tableboxILaplaceTWO}}
 \caption{\label{tab:ILaplaceTWO}  \textsc{I\_Laplace}($n$,2). Comparison of the CPU time and
solution accuracy among CSVD($L=I$), CGSVD and RGSVD, with
sample sizes $l=$ 150, 300,
600 respectively. }
\end{table}

%-----------------------------------------------------------
%-----I_Laplace(n,4)

\newsavebox{\tableboxILaplaceFOUR}
\begin{lrbox}{\tableboxILaplaceFOUR}

\begin{tabular}{cccccccccccc}\hline
& \multicolumn{3}{c}{CSVD} & & \multicolumn{3}{c}{ CGSVD } & & \multicolumn{3}{c}{ RGSVD }\\
 \cline{2-4} \cline{6-8} \cline{10-12}
 $n$  & $T(s)$ & $\mu$ & $err$ & &   $T(s)$ & $\mu$ & $err$  & &  $T(s)$ & $\mu$ & $err$ \\
\hline

                500   & 0.41 & 8.31E-05 & 7.56E-01 & & 0.47  & 6.75E-04 &  4.43E-01  & & 0.20  & 2.44E-01 &  4.94E-02 \\
                1000  & 1.27 & 2.52E-04 & 7.49E-01 & & 2.50  & 2.20E-03 &  3.49E-01  & & 0.99  & 2.56E-01 &  3.90E-02 \\
                2000  & 8.34 & 1.66E-04 & 7.45E-01 & & 19.0  & 3.80E-03 &  2.53E-01  & & 5.41  & 2.50E-01 &  3.22E-02 \\

\hline
\end{tabular}

\end{lrbox}

\begin{table}
 \centering \scalebox{0.8}{\usebox{\tableboxILaplaceFOUR}}
 \caption{\label{tab:ILaplaceFOUR}  \textsc{I\_Laplace}($n$,4). Comparison of the CPU time and
solution accuracy among CSVD($L=I$), CGSVD and RGSVD, with
sample sizes $l=$ 150, 300,
600 respectively. }
\end{table}

{\bf Example 3} (\textsc{Foxgood}). This example arises from discretization of
the integral equation \eqref{Fredholm1st} with the
kernel $K(s,t)=(s^2+t^2)^\frac{1}{2}$ and the exact solution
$x(t)=t$ \cite{Hansen_regutool}.

\begin{figure}[htbp] \centering
\includegraphics[width=\textwidth]{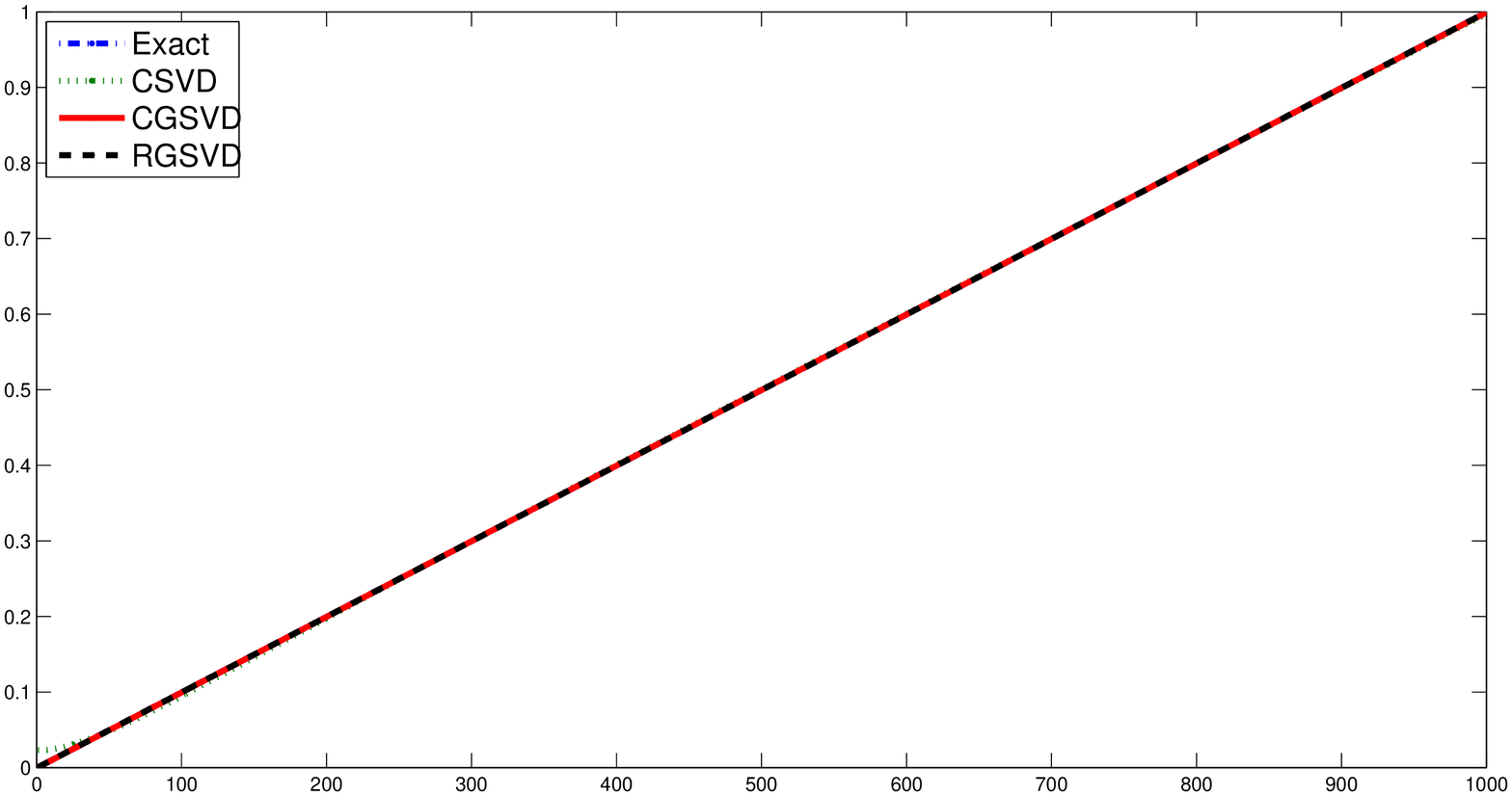}
\caption{  \textsc{Foxgood} matrix of size $n=1000$. }
\label{fig_foxgood}
\end{figure}

The numerical results can be found in Figure \ref{fig_foxgood} and Table \ref{tab:foxgood}.
Similar conclusions and observations can be drawn for this example and
the remaining 3 examples as for Example 1, so we will not repeat them.
The CPU times and the solution accuracy clearly demonstrate the advantages of the new method.

\newsavebox{\tableboxFoxgood}
\begin{lrbox}{\tableboxFoxgood}

\begin{tabular}{cccccccccccc}\hline
& \multicolumn{3}{c}{CSVD} & & \multicolumn{3}{c}{ CGSVD } & & \multicolumn{3}{c}{ RGSVD }\\
 \cline{2-4} \cline{6-8} \cline{10-12}
 $n$  & $T(s)$ & $\mu$ & $err$ & &   $T(s)$ & $\mu$ & $err$  & &  $T(s)$ & $\mu$ & $err$ \\
\hline

                500   & 0.50 & 2.66E-04 & 8.30E-03 & & 0.48  & 7.03E+01 &  3.09E-04  & & 0.07  & 1.41E+02 &  1.25E-04 \\
                1000  & 1.59 & 4.71E-04 & 5.00E-03 & & 2.54  & 2.81E+02 &  1.40E-04  & & 0.10  & 3.90E+05 &  1.73E-05 \\
                2000  & 9.42 & 2.66E-04 & 2.70E-03 & & 18.9  & 1.12E+03 &  1.77E-04  & & 0.21  & 1.35E+03 &  1.45E-04 \\

\hline
\end{tabular}

\end{lrbox}

\begin{table}
 \centering \scalebox{0.8}{\usebox{\tableboxFoxgood}}
 \caption{\label{tab:foxgood}  \textsc{Foxgood}. Comparison of the CPU time and
solution accuracy among CSVD($L=I$), CGSVD and RGSVD.  }
\end{table}

{\bf Example 4} (\textsc{Gravity}). This example arises from discretization of
the equation \eqref{Fredholm1st} with the kernel $K(s,t)=d[
d^2+(s-t)^2]^{-\frac{3}{2}}$,
a one-dimensional model problem in gravity surveying.
Here $d$ is the depth of the point
source and controls the decay of the singular values
\cite{Hansen_regutool}.

\begin{figure}[htbp] \centering
\includegraphics[width=\textwidth]{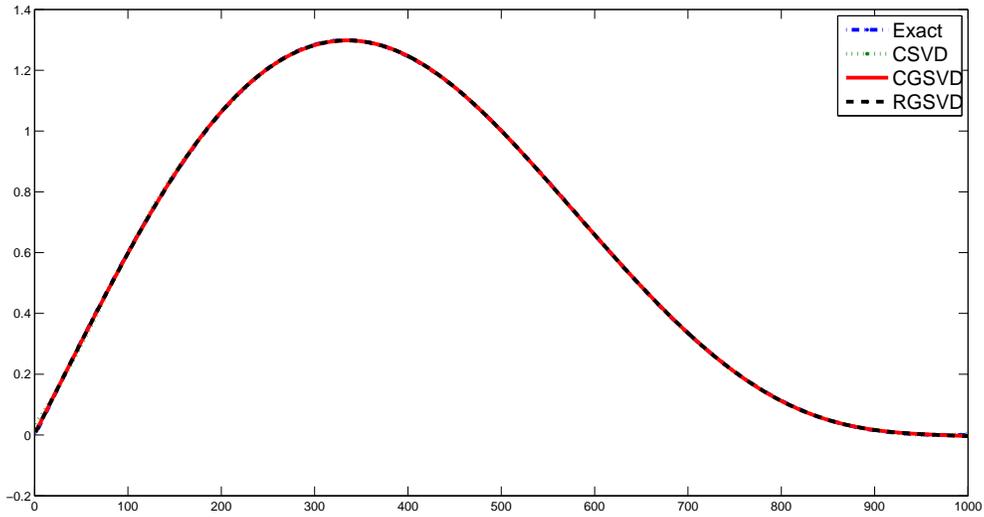}
\caption{  \textsc{Gravity} matrix of size $n=1000$. }
\label{fig_gravity}
\end{figure}

The numerical results are given in Figure \ref{fig_gravity} and Table
\ref{tab:gravity}, when we take $d=0.25$ and the exact solution
$x(t)=\sin(\pi t) +0.5 \sin( 2 \pi t)$.

\newsavebox{\tableboxGravity}
\begin{lrbox}{\tableboxGravity}

\begin{tabular}{cccccccccccc}\hline
& \multicolumn{3}{c}{CSVD} & & \multicolumn{3}{c}{ CGSVD } & & \multicolumn{3}{c}{ RGSVD }\\
 \cline{2-4} \cline{6-8} \cline{10-12}
 $n$  & $T(s)$ & $\mu$ & $err$ & &   $T(s)$ & $\mu$ & $err$  & &  $T(s)$ & $\mu$ & $err$ \\
\hline

                500   & 0.40 & 5.80E-03 & 5.10E-03 & & 0.44  & 9.69E+00 &  1.30E-03  & & 0.07  & 9.69E+00 &  1.30E-03 \\
                1000  & 1.32 & 4.10E-03 & 4.60E-03 & & 2.52  & 1.90E+01 &  7.26E-04  & & 0.09  & 1.90E+01 &  7.26E-04 \\
                2000  & 9.29 & 4.00E-03 & 4.90E-03 & & 18.5  & 6.69E+01 &  1.20E-03  & & 0.19  & 6.69E+01 &  1.20E-03 \\

\hline
\end{tabular}

\end{lrbox}

\begin{table}
 \centering \scalebox{0.8}{\usebox{\tableboxGravity}}
 \caption{\label{tab:gravity}  \textsc{Gravity}. Comparison of the CPU time and
solution accuracy among CSVD($L=I$), CGSVD and RGSVD.  }
\end{table}

\begin{figure}[htbp] \centering
\includegraphics[width=\textwidth]{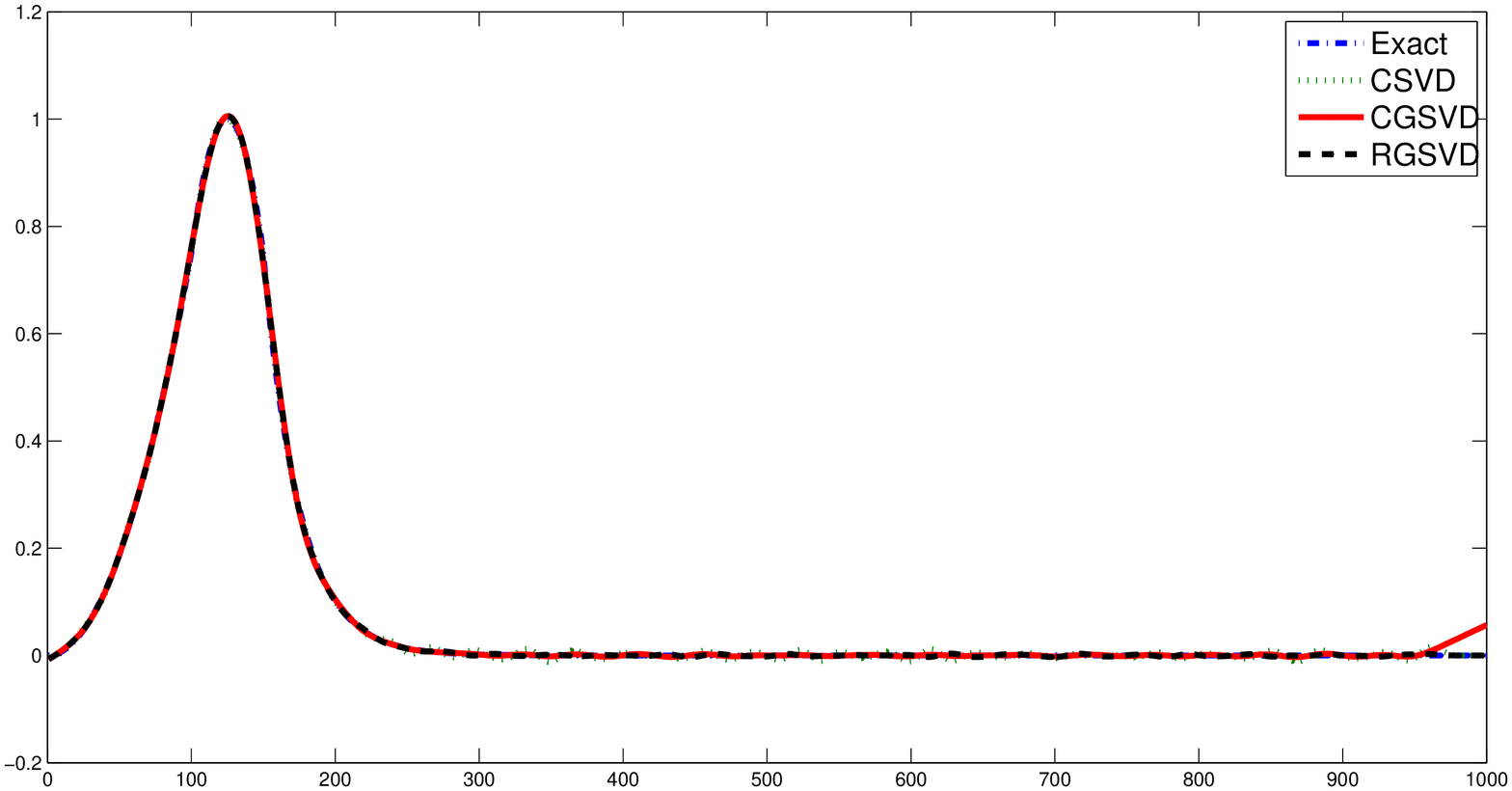}
\caption{  \textsc{Heat} matrix of size $n=1000$. } \label{fig_heat}
\end{figure}

{\bf Example 5} (\textsc{Heat}). This example is the discretization of a Volterra
integral equation of the first kind with the kernel $K(s,t)=k(s-t)$,
where $k(t)=\frac{t^{-2/3}}{2\sqrt{\pi}}exp(-\frac{1}{4t})$;
see \cite{Hansen_regutool} for more detail.
The numerical results are given in Figure \ref{fig_heat} and Table \ref{tab:heat}.

\newsavebox{\tableboxHeat}
\begin{lrbox}{\tableboxHeat}

\begin{tabular}{cccccccccccc}\hline
& \multicolumn{3}{c}{CSVD} & & \multicolumn{3}{c}{ CGSVD } & & \multicolumn{3}{c}{ RGSVD }\\
 \cline{2-4} \cline{6-8} \cline{10-12}
 $n$  & $T(s)$ & $\mu$ & $err$ & &   $T(s)$ & $\mu$ & $err$  & &  $T(s)$ & $\mu$ & $err$ \\
\hline
                500   & 0.46 & 1.41E-04 & 1.62E-02 & & 0.48  & 4.00E-03 &  3.67E-02  & & 0.06  & 2.90E-03 &  1.59E-02 \\
                1000  & 1.79 & 7.97E-05 & 2.23E-02 & & 2.59  & 1.03E-02 &  3.14E-02  & & 0.09  & 1.11E-02 &  1.19E-02 \\
                2000  & 9.53 & 1.17E-04 & 1.23E-02 & & 18.7  & 4.14E-02 &  1.32E-02  & & 0.19  & 1.31E-02 &  1.57E-02 \\

\hline
\end{tabular}

\end{lrbox}

\begin{table}
 \centering \scalebox{0.8}{\usebox{\tableboxHeat}}
 \caption{\label{tab:heat}  \textsc{Heat}. Comparison of the CPU time and
solution accuracy among CSVD($L=I$), CGSVD and RGSVD.  }
\end{table}

\begin{figure}[htbp] \centering
\includegraphics[width=\textwidth]{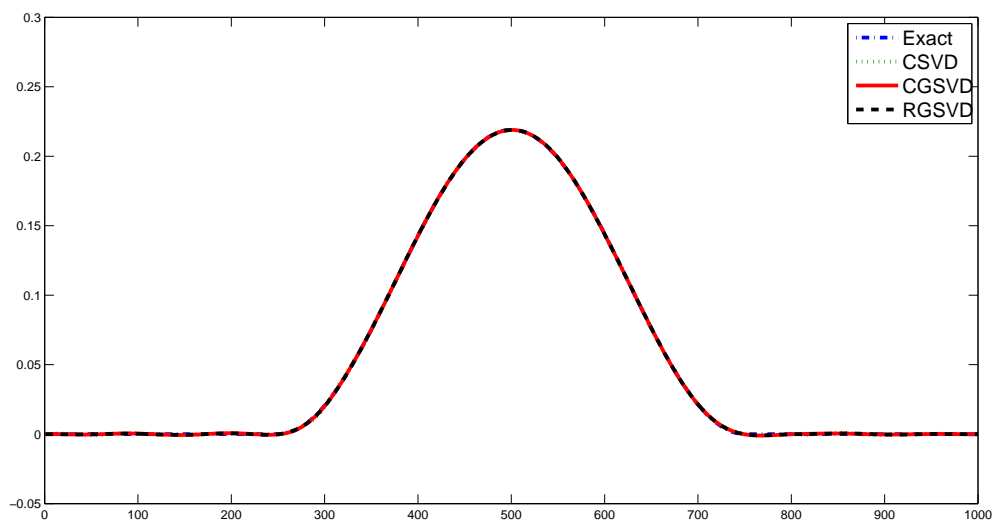}
\caption{  \textsc{Phillips} matrix of size $n=1000$. }
\label{fig_phillips}
\end{figure}

{\bf Example 6} (\textsc{Phillips}). This last testing problem
arises from the discretization of the Fredholm integral equation of
the first kind \eqref{Fredholm1st} designed by D. L. Phillips;
see \cite{Hansen_regutool} for more detail.
The numerical results are given in Figure \ref{fig_phillips} and Table \ref{tab:phillips}.

\newsavebox{\tableboxPhillips}
\begin{lrbox}{\tableboxPhillips}

\begin{tabular}{cccccccccccc}\hline
& \multicolumn{3}{c}{CSVD} & & \multicolumn{3}{c}{ CGSVD } & & \multicolumn{3}{c}{ RGSVD }\\
 \cline{2-4} \cline{6-8} \cline{10-12}
 $n$  & $T(s)$ & $\mu$ & $err$ & &   $T(s)$ & $\mu$ & $err$  & &  $T(s)$ & $\mu$ & $err$ \\
\hline

                500   & 0.51 & 8.70E-03 & 4.00E-03 & & 0.49  & 1.24E-00 &  3.00E-03  & & 0.06  & 1.15E-00 &  2.80E-03 \\
                1000  & 1.59 & 6.40E-03 & 5.40E-03 & & 2.89  & 2.88E-00 &  3.20E-03  & & 0.09  & 2.71E-00 &  2.90E-03 \\
                2000  & 10.7 & 5.30E-03 & 4.40E-03 & & 20.3  & 8.76E-00 &  2.60E-03  & & 0.19  & 7.68E-00 &  2.22E-03 \\

\hline
\end{tabular}

\end{lrbox}

\begin{table}
 \centering \scalebox{0.8}{\usebox{\tableboxPhillips}}
 \caption{\label{tab:phillips}  \textsc{Phillips}. Comparison of the CPU time and
solution accuracy among CSVD($L=I$), CGSVD and RGSVD.  }
\end{table}

\section{Conclusion}

We have considered the randomized algorithms for the solutions
of discrete ill-posed problems in general form. Several strategies are discussed
to transform the problem of general form into the standard one, then the
randomized strategies in \cite{XiangZou_InverseProb13} can be
applied. The second approach we have proposed is to work on the problem of general form
directly. We first reduce the original large-scale problem essentially
by using the randomized algorithm RGSVD, so flops and
memory are significantly saved.
Our numerical experiments show that, using RGSVD we can still achieve
the approximate regularized solutions of the same accuracy as the classical GSVD,
but gain obvious robustness, stability and computational time as we need only to work on
problems of much smaller size.

\end{document}